\documentclass[a4paper,10pt]{article}

\usepackage{mathrsfs}
\usepackage{amsmath}
\usepackage{graphics}
\usepackage{amssymb}
\usepackage{hyperref} 
\date{}
\input xy 
\xyoption{all}  
\input{xypic} 

\newcommand{\maths}[1]{{\rm #1}}

\renewcommand{\leq}{\leqslant}
\renewcommand{\geq}{\geqslant}

\newcommand{\sch}[1]{\mathscr #1}

\renewcommand{\hat}{\widehat}

\renewcommand{\phi}{\varphi}

\newcommand{\got}[1]{{\mathfrak #1}}
\renewcommand{\Bbb}{\mathbb}

\newcommand{\RR}{{\Bbb R}}

\newcommand{\PP}{{\Bbb P}}
\newcommand{\Aff}{{\Bbb A}}

\newcommand{\NN}{{\Bbb N}}

\newcommand{\QQ}{{\Bbb Q}}

\renewcommand{\epsilon}{\varepsilon}

\newcommand{\an}{^{\rm an}}

\newcommand{\zero}{^{\mbox{\tiny o}}}
\newcommand{\zeroo}{^{\mbox{\tiny oo}}}

\newcommand{\spec}{{\rm Spec}\;}

\newcommand{\red}{\widetilde}

\newcommand{\hres}{{\sch H}}

\newcommand{\deux}[1]{\refstepcounter{subsection}\label{#1}\medskip\noindent {\bf (\thesubsection)}\hspace{.1cm}}
\newcommand{\trois}[1]{\refstepcounter{subsubsection}\label{#1}\medskip\noindent {\bf
    (\thesubsubsection)}\hspace{.1cm}}   
   \author{\sc Antoine Ducros\thanks{During this work the author was partially supported by ValCoMo (ANR-13-BS01-0006), and by a Rosi and Max Varon visiting professorship at Weizmann Institute.}\\ \small IMJ-PRG \\ \small Institut universitaire de France}
   \title{About Hrushovski and Loeser's work on the homotopy type of Berkovich spaces}

   \date{}
   
\begin{document}
   \maketitle
   
\setcounter{section}{-1}
\section{Introduction}

\deux{intro-berko}
At the end of the eighties, V. Berkovich developed in \cite{brk1}
and \cite{brk2} a new theory of analytic geometry over non-archimedean fields, coming
after those by Krasner, Tate \cite{tate}  and Raynaud \cite{ray}. 
One of the main advantages of this approach is that
the resulting spaces enjoy very nice {\em topological}
properties: they are locally compact and locally {\em pathwise} connected (although non-archimedean fields
are totally disconnected, and often not locally compact).
Moreover, Berkovich spaces have turned out to be ``tame" objects -- in the informal sense of Grothendieck's {\em Esquisse d'un programme}. Before
illustrating
this rather vague assertion
by several examples, let us fix some terminology
and notations.

\deux{terminology}
Let~$k$ be a field which is complete with respect to a non-archimedean absolute value $|.|$ and let
$k\zero$ be the valuation ring $\{x\in k, |x|\leq 1\}$. 

\trois{anal-space}
A (Berkovich) analytic space over $k$ is a (locally compact, locally path-connected)
topological space $X$ 
equipped with some extra-data, among which: 

- a sheaf of $k$-algebras whose sections are called {\em analytic functions}, which makes $X$ a locally ringed space; 

- for every point $x$ of $X$, a complete non-archimedean field $(\hres(x),|.|)$ endowed with 
an isometric embedding $k\hookrightarrow \hres(x)$, and an "evaluation" morphism
from $\sch O_{X,x}$ to $\hres(x)$, which is denoted by $f\mapsto f(x)$. 

\trois{analytify}To every $k$-scheme of finite type $X$ is associated in a functorial way
a $k$-analytic space $X\an$, which is called its {\em analytification} (\cite{brk2}, \S 2.6); it is provided
with a natural morphism of locally ringed spaces $X\an \to X$
 If $f$ is a section of $\sch O_X$, we will
often still denote by $f$ its pull-back to $\sch O_{X\an}(X\an)$.

\medskip
Let us recall what the underlying topological space of $X\an$ is. As a set, it consists
of couples $(\xi,|.|)$ where $\xi \in X$ and where $|.|$ is a non-archimedean 
absolute value on $\kappa(\xi)$ extending the given absolute value of $k$. If $x=(\xi,|.|)$ is a point
of $X\an$, then $\hres(x)$ is the completion of the valued field $(\kappa(\xi),|.|))$.
The map $X\an \to X$
is nothing but $(\xi,|.|)\mapsto \xi$; note that the pre-image of a Zariski-open subset $U$
of $X$ on $X\an$ can be identified with $U\an$.  The topology $X\an$ is equipped with is the coarsest such than:

\medskip
$\bullet$ $X\an \to X$ is continuous (otherwise said, $U\an$ is an open subset of $X\an$ for every Zariski-open subset $U$ of $X$); 

$\bullet$ for every Zariski-open subset $U$ of $X$ and every $f\in \sch O_X(U)$, the 
map $(\xi,|.|)\mapsto |f(\xi)|$ from $U\an$ to $\RR_+$ is continuous.

\medskip
Let us now assume that $X$ is affine, say $X=\spec A$. The topological space underlying $X\an$ can then
be given another description: it is
the set of all multiplicative maps $\phi \colon A\to \RR_+$ that extend the absolute value of $k$
and that satisfy the inequality $\phi(a+b)\leq \max(\phi(a), \phi(b))$ for every $(a,b)\in A^2$ (such a map
will be simply called a {\em multiplicative semi-norm}); its topology is the one inherited
from the product topology on $\RR_+^A$. 
This is related as follows to the former description: 

\medskip
$\bullet$ to any couple $(\xi,|.|)$ as above corresponds 
the multiplicative semi-norm $\phi$ defined by the formula $\phi(f)= | f(\xi)|$; 

$\bullet$ to any multiplicative semi-norm $\phi$ corresponds the couple $(\got p, |.|)$ where $\got p$
is the kernel of $\phi$ (this is a prime ideal of $A$) and $|.|$ is the absolute value on ${\rm Frac}\;A/\got p$
extending the multiplicative norm on $A/\got p$ induced by $\phi$. 

\medskip
Modulo this new description the map $X\an \to X$ simply sends a semi-norm to its kernel.

\trois{semi-alg}{\bf Definition}. 
Let $X$ be a $k$-scheme of finite type. A subset of $X\an$ is said to 
be {\em semi-algebraic}
if it can be defined, locally for the Zariski-topology of~$X$,
by a boolean combination of inequalities
$|f|\Join \lambda |g|$ where~$f$ and~$g$ are 
sections of $\sch O_X$, where~$\Join\in\{<,>,\leq,\geq\}$,
and where $\lambda\in \RR_+$; it is called {\em strictly}
semi-algebraic if the~$\lambda$'s can be chosen in~$|k|$.

\trois{formal-scheme}
To every $k\zero$-formal scheme $\got X$ (topologically) of 
finite type is associated in a functorial way a compact $k$-analytic space
$\got X_\eta$, which is called its {\em generic fiber}
(\cite{vc}, \S 1). 

In what follows, we will be interested in the generic 
fiber $\got X_\eta$ for $\got X$ being a {\em polystable}
$k\zero$-formal scheme. This notion has been introduced by Berkovich
(def. 1.2 of \cite{sm1}); let us simply emphasize here
that semi-stable
$k\zero$-formal schemes are polystable. 

\trois{polyhedron}
{\bf Definition}. A {\em polyhedron}
is a subset of $\RR^n$ (for some $n$) 
which is a finite union of {\em rational}
simplices. 

\deux{tameness}
{\bf Examples of
tameness properties of Berkovich spaces}.

\trois{genfib-polystable}
Let $\mathfrak X$ be a polystable $k\zero$-formal scheme. 
Berkovich proves in~\cite{sm1}
that its generic fiber $\mathfrak X_\eta$
admits a strong deformation retraction to 
one of its closed subset $S(\mathfrak X)$
which is homeomorphic to a polyhedron\footnote{It is more precisely
homeomorphic to the (abstract) ``incidence polyhedron" of the special fiber $\mathfrak X_s$
that encodes the combinatorics of its singularities, {\em e.g}
it is reduced to a point if $\mathfrak X_s$ is smooth and irreducible.}
of dimension $\leq \dim \got X$.

\trois{sm-loc-cont}
Smooth $k$-analytic spaces are locally contractible; this is also proved in~\cite{sm1} by Berkovich, by reduction to
assertion~\ref{genfib-polystable}
above through de Jong's alterations.

\trois{semi-alg}
Let~$X$ be a $k$-scheme of finite type. Every semi-algebraic subset of~$X^{\rm an}$ 
has finitely many connected components, each of which is
semi-algebraic; this was proved by the author in~\cite{semialg}.

\trois{variation-pi0}
Let~$X$ be a compact $k$-analytic space and let~$f$ be an analytic function on~$X$; for every~$\epsilon\geq 0$, denote by~$X_{\epsilon}$ the
set of~$x\in X$ such that~$|f(x)|\geq \epsilon$. There exists a finite partition~$\sch P$ of~$\RR_+$ in intervals
such that for every~$I\in \sch P$ and every~$(\epsilon,\epsilon')\in I^2$ with~$\epsilon\leq \epsilon'$,
the natural map~$\pi_0(X_{\epsilon'})\to \pi_0(X_\epsilon)$
is bijective. This has been established by Poineau in~\cite{poin} (it 
had already been proved in the particular case where~$f$ is invertible by Abbes and Saito in~\cite{abb}).

\deux{presentation-hl}
{\bf Hrushovski and Loeser's work}.
As far as tameness
is concerned, there has been in 2009-2010
a major breakthrough, namely the work~\cite{hl},
by Hrushovski and Loeser. 
Let us quickly explain what their main results consist of. We fix
a {\em quasi-projective}
$k$-scheme and a semi-algebraic subset $V$ of $X\an$. 
Hrushovski and Loeser 
have proven the following.

\trois{def-ret-hl}
There exists a strong deformation retraction from $V$ to one of its closed subsets which is homeomorphic
to a polyhedron. More precisley, Hrushovski and Loeser build a continuous map~$h: [0;1]\times V \to V$,
and prove the existence of a compact subset $S$ of~$V$ homeomorphic to 
a polyhedron such that the following hold:

\begin{itemize}

\item[i)] $h(0,v)=v$ and~$h(1,v)\in S$ for all~$v\in V$;

\item[ii)] $h(t,v)=v$ for all~$t\in [0;1]$ and all~$v\in S$; 

\item[iii)]  $h(1, h(t,v))=h(1,v)$ for all~$t\in [0;1]$ and all~$v\in V$. 

\end{itemize} 

\trois{loc-cont-hl}
The topological space~$V$ is locally contractible. 

\trois{homtypefib-hl}
Let $Y$ be a $k$-scheme of finite type and let $\phi \colon X\to Y$ be a $k$-morphism. The set of homotopy types 
of fibers of the map~$\phi\an|_V : V\to Y\an$ is finite\footnote{One could ask the same
question concerning the set of {\em homeomorphism} types; to the author's knowledge, the answer is
not known, and likely not easy to obtain using Hrushovski and Loeser's approach.}.

\trois{functionslevel-hl}Let~$f$ be a function belonging to~$\sch O_X(X)$. For every~$\epsilon\geq 0$ let us denote by~$V_\epsilon$ the set of~$x\in V$
such that~$|f(x)|\geq \epsilon$.  There exists a finite partition~$\sch P$ of~$\RR_+$ in intervals
such that for every~$I\in \sch P$ and every~$(\epsilon,\epsilon')\in I^2$ with~$\epsilon\leq \epsilon'$,
the embedding~$V_{\epsilon'}\hookrightarrow V_\epsilon$
is a homotopy equivalence.

\deux{commentaires-hl} {\bf Comments.} 

\trois{loc-contract-corollaire}
Since~$V$ has a basis of open subsets
which are semi-algebraic subsets of~$X\an$, 
assertion~\ref{loc-cont-hl}
is an obvious corollary of assertion~\ref{def-ret-hl}
and of the local 
contractibility of polyhedra; note that the property iii)
of~\ref{loc-cont-hl}
is here crucial. 

\trois{smooth-previous}
Even when~$X$ is smooth, projective and when~$V=X\an$, assertion~\ref{def-ret-hl}
was 
previously known only when~$X$ has 
a polystable model. 

\trois{qproj-why}
The quasi-projectivity assumption on~$X$
is needed by Hrushovski and Loeser for {\em technical}
reasons, but everything remains likely true for $X$ any $k$-scheme of finite type.  

\trois{content-hl}
In order to prove these tameness result, 
Hrushovski and Loeser develop a new kind of geometry over valued fields
of arbitrary height, using
highly sophisticated tools from model theory. 
In fact, most of their work
is actually devoted to this geometry, and they transfer their results to the Berkovich
framework only at the very end of their paper (chapter 14 of \cite{hl});
for instance, they
first prove counterparts of assertions \ref{def-ret-hl} - \ref{functionslevel-hl}
in their setting. 

\deux{purpose}
{\bf About this text}. 
In this article we survey the methods and results of~\cite{hl}\footnote{The reader may also refer
to the more detailed survey \cite{bourb}.}. We will describe roughly
the geometry of Hrushovski and Loeser, 
say a few words about the links
between their spaces and those of Berkovich, 
and 
give a coarse sketch of the proof
of their version of assertion~\ref{def-ret-hl}. 

In the last section, we will explain how a key finiteness result of~\cite{hl} 
has been used by the author in~\cite{polyt}, to show the following: 
if~$X$ is an~$n$-dimensional analytic space and if~$f: X\to {\mathbb G}_{m,k}\an$ is any morphism, 
the pre-image of the ``skeleton'' of~${\mathbb G}_{m,k}^{n,\rm an}$
under~$f$
inherits a canonical 
piecewise-linear structure. 

\subsection*{Acknowledgements}
I'm grateful to the referees for their careful reading, and for many helpful
comments and suggestions. 

\section{Model theory of valued fields : basic definitions}

\deux{notation-valued fields}
{\bf General conventions about valued fields}. 

\trois{krull-valuation}
In this text, a {\em valuation}
will be an abstract Krull valuation, neither assumed
to be of height one nor non-trivial.  
We will use the {\em multiplicative notation}: if $k$ is a valued field then unless otherwise stated
its valuation will be denoted by $|.|$, and
its value group will not be given a specific name -- we will simply use the notation $|k^\times|$, 
and $1$ for its unit element. 
We will formally add to 
the ordered group $|k^\times|$ an absorbing, smallest element $0$ and set $|0|=0$; with these conventions
one then has $|k|=\{0\}\cup |k^\times|$ and 
$$|a+b|\leq \max (|a|, |b|)$$ for every $(a,b)\in k^2$.

We will denote by $|k^\times|^{\QQ}$ the divisible hull of the ordered group $|k^\times|$, and
set for short $|k|^{\QQ}=\{0\}\cup |k^\times|^{\QQ}$.

\medskip
The reader should be aware that the multiplicative notation, which is usual in Berkovich's theory, 
is {\em not}
consistent with Hrushovski and Loeser's choice. Indeed, they use the additive notation, 
call ``$0$" what we call ``$1$", call ``$+\infty$" what we call ``$0$", and their order
is the opposite of ours. 

\trois{residue-field}
If $k$ is a valued field, we will set 
$$k\zero=\{\in k, |x|\leq 1\}\;{\rm and}\;k\zeroo=\{x\in k, |x|<1\}.$$ Then $k\zero$ is a subring of $k$; it is local,
and $k\zeroo$ is its unique maximal ideal. The quotient $k\zero/k\zeroo$ will be denoted by $\tilde k$; it is called the {\em residue
field}
of $(k,|.|)$ (or simply of $k$ if there is no ambiguity about the valuation).

\deux{deffun}
{\bf Definable functors: the embedded case}. 
Let~$k$ be a valued field. 
We fix an algebraic closure~$\bar k$ of~$k$, and an extension of the valuation
of~$k$ to~$\bar k$. Let~$\mathsf M$ 
be the category of algebraically closed valued extensions of~$\bar k$
whose valuation is non-trivial (morphisms are isometric~$\bar k$-embeddings). 

\trois{def-subset}
{\bf Definable \em subsets}. 
Let~$X$ and~$Y$ be~$k$-schemes of finite type and
let~$F\in \mathsf M$. 
A subset of~$X(F)$ will be said to be {\em $k$-definable}
if it can be defined, locally for the Zariski-topology of~$X$, by a boolean combination of
inequalities of the form
$$|f|\Join\lambda |g|$$ where~$f$ and~$g$ are regular functions,
where~$\lambda \in |k|$, and where~$\Join\in\{\leq,\geq,<,>\}$. A map
from a~$k$-definable subset of~$X(F)$ to a~$k$-definable subset of~$Y(F)$
will be said to be~$k$-definable if its graph is
a~$k$-definable subset of~$X(F)\times Y(F)=(X\times_k Y)(F)$. 

\trois{comment-definition}
{\em Comments}. Definability is in fact a general notion of model theory
which 
makes sense with respect to any given {\em language}. The one we have introduced here is
actually definability in the {\em language of valued fields}, and
the reader is perhaps more familiar with the 
two following examples: 
definability in the language of fields, which is nothing but Zariski-constructibility; 
and definability in the language of ordered fields, which is nothing but semi-algebraicity.

\trois{def-subfunctor}
{\bf Definable \em sub-functors}.
The scheme~$X$ induces a
functor~$F\mapsto X(F)$
from~$\mathsf M$ to $\mathsf{Sets}$, which we will also
denote by~$X$. A sub-functor~$D$ of~$X$ will be said to be {\em $k$-definable}
if it can be defined, locally for the Zariski-topology of~$X$, by a boolean combination of
inequalities of the form
$$|f|\Join\lambda |g|$$ where~$f$ and~$g$ are regular functions,
where~$\lambda \in |k|$, and where~$\Join\in\{\leq,\geq,<,>\}$. A natural transformation
from a~$k$-definable sub-functor of~$X$ to a~$k$-definable sub-functor~of
$Y$
will be said to be~$k$-definable if its graph is a~$k$-definable sub-functor of~$X\times_k Y$. 
If~$D$ is a~$k$-definable sub-functor of~$X$, then~$D(F)$ is for every~$F\in \mathsf M$
a~$k$-definable subset of~$X(F)$; if~$D'$ is a~$k$-definable sub-functor of~$Y$ and if~$f: D\to D'$
is a~$k$-definable natural transformation, $f(F): D(F)\to D'(F)$
is for every~$F\in \mathsf M$ a~$k$-definable map.

The following fundamental facts are straightforward consequences 
of the so-called
{\em quantifier elimination for non-trivially valued algebraically
closed fields}. It is usually attributed to Robinson, though the result 
proved in his book \cite{robi}
is weaker -- but the main
ideas are there; for a complete proof, {\em cf.}
\cite{zoe}, \cite{pres}
or~\cite{weis}.

\medskip
i) For every~$F\in \mathsf M$, the assignment~$D\mapsto D(F)$ establishes
a one-to-one correspondence between $k$-definable {\em sub-functors} of~$X$ and
$k$-definable {\em subsets}
of~$X(F)$. The analogous result holds for $k$-definable natural transformations
and~$k$-definable maps. 

ii) If~$f: D\to D'$ is a $k$-definable natural transformation between two~$k$-definable
sub-functors~$D\subset X$ and~$D'\subset Y$,
the sub-functor of~$D'$ that sends~$F$ to the 
image of~$f(F): D(F)\to D'(F)$ is~$k$-definable. 

\medskip
{\em Comments about assertion i).}
For a given~$F\in \mathsf M$, 
assertion~i) allows to identify~$k$-definable subsets of~$X(F)$
with~$k$-definable
sub-functors of~$X$, and the choice of one of those viewpoints
can be somehow a matter of taste. But even if one is
actually only interested
in~$F$-points, 
it can be useful to think
of a~$k$-definable subset~$D$
of~$X(F)$
as a sub-functor of~$X$. Indeed, this allows 
to consider points of~$D$ 
over valued fields 
larger than~$F$, which often encode in a natural
and efficient way the ``limit'' behavior of~$D$. 

\trois{subvar}
{\em Subvarieties of $X$ and definable subfunctors}.
Here we consider only {\em field-valued}
points; for that reason, 
this theory is totally insensitive to nilpotents.
For example, $X_{\rm red}$ and $X$ define the same functor in our setting.
But let us now give a more subtle example of this 
phenomenon. 
If~$Y\hookrightarrow X$ is a $k$-subscheme, 
then $F\mapsto Y(F)$ is obviously a $k$-definable subfunctor of $X$. 

Now, let $L$ be the perfect closure of $k$
(which embeds canonically in any $F\in \mathsf M$), and let $Y$ 
be an $L$-subscheme of $X_L$. The functor 
$F\mapsto Y(F)$ is then again $k$-definable, even if $Y$ can not
be defined over $k$ in the usual, algebro-geometric sense. Indeed, 
to see it we immediately reduce to the case where $X$ is affine ; now, 
$Y$ is defined by an ideal $(f_i)_{1\leq i\leq n}$, where the $f_i$'s all
belong to $\mathscr O_{X_L}(X_L)$.  But then, again since we only 
consider field-valued points, the subfunctor $F\mapsto Y(F)$ is equal 
to the one defined by the ideal $(f_i^{p^n})$ for any $n$, where $p$ is the
characteristic exponent of $k$. By taking $n$ large enough so that
$f_i^{p^n}\in \mathscr O_X(X)$ for every $i$, we see that our functor 
is $k$-definable. 

\deux{abstract-defin}
{\bf Definable functors: the abstract case.}

\trois{abstract-embedable}
{\bf Definably embeddable functors}. Let us say
that a
functor~$D$ from~$\mathsf M$ to~$\mathsf{Sets}$ is
{\em $k$-definably embeddable} 
if there exists a $k$-scheme of finite type~$X$,
a definable sub-functor~$D_0$ of~$X$ and an isomorphism~$D\simeq D_0$. 

The functor~$D_0$ is
not 
{\em a priori}
uniquely determined
(up to unique $k$-definable isomorphism). But for all functors~$D$ we will consider below, 
we will be
given
{\em implicitly}
not only~$D(F)$ for every~$F\in \mathsf M$, but also $D(T)$ for every~$F\in \mathsf M$
and every~$F$-definable sub-functor~$T$ of an~$F$-scheme of finite type;
otherwise said, if~$D$ classifies 
objects of a certain kind, we not only know what such an object
defined over~$F$ is, but also what an embedded~$F$-definable family of such objects is. 
It follows from Yoneda's lemma that such data ensure the canonicity of~$D_0$ as soon 
as it exists (see~\cite{bourb}, \S 1 for more detailed explanations about those issues).

\trois{limite-definable}
It turns out that the class of $k$-definably embeddable functors is 
too restrictive, 
for the following reason. Let $D$ be such a functor and let $R$ be a subfunctor of $D\times D$. Assume
that $R$ is itself $k$-definably embeddable, and that it is an equivalence relation -- that is, that $R(F)\subset D(F)\times D(F)$
is the graph of an equivalence relation on $D(F)$ for every $F\in \mathsf M$. 
The quotient functor
$$D/R:=F\mapsto D(F)/R(F)$$ is then not necessarily $k$-definably embeddable. In order to remedy this, 
one simply enlarges our class of functors by {\em forcing}
it to be stable under such quotient operations. 

\trois{def-definable}
{\bf Definition}. A functor $\Delta \colon \mathsf M\to \mathsf{Sets}$ is said to be {\em $k$-definable}
if it is isomorphic to $D/R$ for some $k$-definably embeddable functor $D$ and some $k$-definably embeddable
equivalence relation $R\subset D\times D$.

A natural transformation between $k$-definable functors is called $k$-definable if its graph
is itself a $k$-definable functor. 

\trois{tautology-quotient}
{\em Remark}. Let $D$ be $k$-definable functor and let $R\subset D\times D$ be a $k$-definable
equivalence relation. It follows from the definition that the quotient functor $D/R$ is $k$-definable.

\trois{comments-abstract-def}
{\em Comments}. Quantifier elimination ensures that there is no conflict of terminology: if $X$ is a $k$-scheme
of finite type, then a subfunctor $D$ of $X$ is $k$-definable in the above sense if and only if it is $k$-definable
in the sense of~\ref{def-subfunctor}.

\medskip
The following counterparts of assertions i) and ii) of {\em loc.~cit.} hold 
(the first one is again a consequence of quantifier elimination; the second one comes from remark~\ref{tautology-quotient}
above). 

$\bullet$ {\em Counterpart of~i)}. Let~$\Delta$ be a~$k$-definable functor. For every~$F\in \mathsf M$, let us say
that a subset of~$\Delta(F)$ is~$k$-definable if it can be written~$D(F)$ for~$D$ a~$k$-definable sub-functor of~$\Delta$. The assignment
~$D\mapsto D(F)$ then induces a one-to-one correspondence between~$k$-definable subsets of~$\Delta(F)$ and~$k$-definable sub-functors
of~$\Delta$. The analogous statement for~$k$-definable transformations holds. 

$\bullet$ {\em Counterpart of~ii).} If~$f: D\to D'$ is a $k$-definable natural transformation between two~$k$-definable
functors
the sub-functor of~$D'$ that sends~$F$ to the 
image of~$f(F): D(F)\to D'(F)$ is~$k$-definable.

\trois{k-def-examples}
{\bf Examples.} The following functors are $k$-definable; their descriptions
as quotients are left to the readers.  

\medskip
\begin{itemize} 
\item[$\bullet$] The functor~$\Gamma : F\mapsto |F^\times|$.

\item[$\bullet$] The functor~$\Gamma_0 : F\mapsto |F|$.

\item[$\bullet$] The functor~$F\mapsto \red F$. 

\item[$\bullet$] If~$a$ and~$b$ are elements of~$|k|^{\QQ}$
with $a\leq b$ the functor~$[a;b]$ 
 that
sends~$F$ to
$$\{c\in |F|,a\leq c\leq b\}$$
is $k$-definable (note that $|k|^{\QQ}\subset |F|$ for every $F\in \mathsf M$).
One defines in an analogous way the functors $[a;b)$, $(a;b)$, $[a;+\infty)$, etc. 
Such functors
are called
{\em $k$-definable intervals}; a $k$-definable interval of the form $[a;b]$ will
be called a $k$-definable {\em segment}. 

\end{itemize} 

\medskip
{\em Remark.}
One can prove that none of the functors
above is $k$-definably embeddable, except the singleton and the empty set (which can be described
as generalized intervals).
What it means can be roughly 
rephrased as follows, say for~$\Gamma_0$ (there is an analogous
formulation for every of the other functors): 
one can not find an algebraic variety~$X$ over~$k$ and 
a natural
way to embed~$|F|$ in~$X(F)$ for every~$F\in \mathsf M$.

\trois{example-definable}
{\bf Example of $k$-definable natural transformations}. Let $I$ 
and $J$ be two $k$-definable intervals. A {\em $k$-monomial map}
from $I$ to $J$ is a natural transformation 
of the form
$$x\mapsto ax^r$$ with $a\in |k|^{\QQ}$ and $r\in \QQ_+$ (resp. $\QQ$)
if $0\in I$ (resp. if $0\notin I$); here we use the convention $0^0=1$. Such a natural 
transformation is $k$-definable. Moreover, for any $k$-definable natural transformation 
$u\colon I\to J$, there exists a finite partition $I=\coprod_n I_n$ of $I$ in $k$-definable
intervals such that $u|_{I_n}$ is $k$-monomial for every $n$.

Assume now that we are given a $k$-definable {\em injection} $u\colon I\to J$, that $I$ contains
an interval of the form $]0;a[$ with $a\in |k^\times|^{\QQ}$ and that $J$ is bounded, {\em i.e}
contained in $[0;R]$ for some $R\in|k^\times|^{\QQ}$. Then we can always decrease $a$ so that $u|_{]0;a[}$ is
$k$-monomial. Since $u$ is injective, the corresponding exponent will be non-zero; since $J$ is bounded, 
it will even be positive. Hence the map $u$ induces an increasing $k$-monomial
bijection between $]0;a[$ and a $k$-definable interval $]0;b[\subset J$.

\deux{more-involved-examples}
{\bf More involved examples of $k$-definable functors}. 

\trois{concatenation} One can
{\em concatenate} 
a finite sequence of $k$-definable segments: one makes the quotient of their disjoint union
by the identification of successive
endpoints and origins; one gets that way
$k$-definable {\em generalized} segments. 

\medskip
The  concatenation of finitely many $k$-definable 
segments {\em with non-zero origin and endpoints}
is $k$-definably isomorphic to a single $k$-definable segment; the proof
consists in writing down an explicit isomorphism and is left to the reader. 

But be aware that this does not
hold in general without our assumption on the origins
and endpoints. For instance, the concatenation $I$ of two copies of~$[0;1]$, where the endpoint $1$
of the first one is identified with the origin $0$ of the second one, is not~$k$-definably
isomorphic 
to a~$k$-definable segment; indeed, it even 
follows
even from~\ref{example-definable}
that there is no $k$-definable
injection from $I$ to a bounded $k$-definable interval.

\trois{definable-pl}
A sub-functor of~$\Gamma_0^n$ is~$k$-definable
if and only if it can be defined by conjunctions and disjunctions
of inequalities of the form 
$$ax_1^{e_1}\ldots x_n^{e_n}\Join bx_1^{f_1}\ldots x_n^{f_n}$$ where $a$ and $b$
belong to $|k|$, and where the $e_i $'s and the $f_i$'s belong to $\NN$, and where $\Join\in\{<,\leq, >, \geq\}$.

\medskip
Any $k$-definable generalized segment is $k$-definably isomorphic to some 
$k$-definable sub-functor of $\Gamma_0^n$
(the proof is left to the reader). 

\trois{ball-definable} {\em The functor~$B$
that sends a field $F\in \mathsf M$ to the set of its closed balls
is~$k$-definable.} Indeed, 
let us
denote by~$T$ 
the~$k$-definable functor 
$$F\mapsto \left\{\left(\begin{array}{cc}
a&b\\0&a\end{array}\right)\right\}_{a\in F^*, b\in F},$$
and by~$T'$ its~$k$-definable
sub-functor $F\mapsto T(F)\cap {\rm GL_2}(F\zero)$.  
The quotient functor~$T/T'$ is~$k$-definable, and
the reader will check that the map that sends $$\left(\begin{array}{cc}
a&b\\0&a\end{array}\right)$$ to the closed ball with center~$b$ 
and radius~$|a|$ 
induces a functorial bijection between~$T(F)/T'(F)$ 
and the set of closed balls of~$F$ of {\em positive} radius. The functor~$B$
is thus isomorphic to 
$F\mapsto \left(T(F)/T'(F)\right) \coprod F$, and is therefore
$k$-definable.

\deux{class-all-defi}
{\bf About the classification of {\em all}
$k$-definable functors.}
By definition,
the description of a $k$-definable functor
involves a $k$-definably embeddable functor
and a $k$-definably embeddable equivalence relation on it. 

\medskip
In fact, in the seminal work
\cite{imaginaries}
Haskell, Hrushovski, and MacPherson prove that the equivalence relation above can always be chosen to be in a particular explicit list, as we now explain.

\trois{def-sn-tn}
Let $n\in \NN$. We denote by $S_n$
the functor that sends $F\in \mathsf M$ to $\mathrm{GL}_n(F)/\mathrm{GL}_n(F^{\circ})$ ; note that $S_1\simeq \Gamma$. 

We denote by $T_n$ be the functor that sends $F\in \mathsf M$
to the quotient of $\mathrm{GL}_n(F)$ by the kernel of
$\mathrm{GL}_n(F^\circ)\to \mathrm{GL}_n(F^\circ/F^{\circ \circ})$. 

By construction, $S_n$ and $T_n$ appear as quotients of $\mathrm{GL}_n$ by $k$-definable 
relations, hence are $k$-definable. 

\trois{desc-defi-general}
Let $X$ be a $k$-scheme of finite type, and let $(n_i)$ and $(m_j)$ be two finite
families of integers. Let $D$ be a $k$-definable sub-functor
of the product $X\times \prod_i \mathrm{GL}_{n_i}\times \prod_j \mathrm{GL}_{m_j}$
The image $\Delta$ of $D$ under
the map $$X\times \prod_i \mathrm{GL}_{n_i}\times \prod_j \mathrm{GL}_{m_j}\to
X\times \prod_i S_{n_i}\times \prod_j T_{m_j}$$ is a $k$-definable sub-functor of
$X\times \prod_i S_{n_i}\times \prod_j T_{m_j}$. Haskell, Hrushovski and Macpherson have
then shown that every $k$-definable functor is isomorphic to such a $\Delta$.

\medskip
\deux{scalar-change} {\bf Ground field extension}. Let $F\in \mathsf M$, let $F_0$ 
be a subfield 
of $F$ containing $k$, and let $\overline{F_0}$ be the algebraic closure
of $F_0$ inside $F$. Let $\mathsf N$ be the category of all non-trivially valued, 
algebraically closed fields containing $\overline{F_0}$; note that $F\in \mathsf N$, 
that $\mathsf N\subset \mathsf M$
(since $\bar k\subset \overline{F_0}$), and that $\mathsf N=\mathsf M$ as 
soon as $F_0$ is algebraic over $k$ (since then $\overline{F_0}=\bar k$). 

\medskip
The notion 
of an $F_0$-definable functor from $\mathsf N$ to $\mathsf {Sets}$, 
and that of an $F_0$-definable subset of $\Delta(L)$ for such a functor $\Delta$
and for $L\in \mathsf N$, make sense. If $D$ is a $k$-definable functor
from $\mathsf M\to \mathsf{Sets}$, its restriction to $\mathsf N$ is
$F_0$-definable ; in particular, the notion of an $F_0$-definable subset
of $D(L)$ for any $L\in \mathsf N$ makes sense.

\section{Hrushovski and Loeser's fundamental construction}

For every~$F\in \mathsf M$, we denote by~$\mathsf M_F$ the category of valued extensions
of~$F$ that belong to~$\mathsf M$. 

\deux{type}
{\bf The notion of a type ({\em cf.} \cite{hl}, \S 2.3-2.9).}
Let~$F\in \mathsf M$ and let~$D$ be an~$F$-definable functor. 

\trois{equiv-type}
Let~$L$ and~$L'$
be two valued fields belonging to~$\mathsf M_F$. 
Let us say that 
a point~$x$ of~$D(L)$ and a point~$x'$  of~$D(L')$ are {\em $F$-equivalent}
if for every $F$-definable sub-functor~$\Delta$ of~$D$, one has
$$x\in \Delta(L)\iff x'\in \Delta(L').$$
Roughly speaking, $x$ and $x'$ are equivalent if they satisfy
the same formulas with parameters in~$F$.

\medskip
We
denote by~$S(D)$ the set of couples~$(L,x)$ where~$L$
belongs to~$\mathsf M_F$
and where~$x$ belongs to~$D(L)$, up to~$F$-equivalence; an element of~$S(D)$
is called a
{\em type}
on~$D$.

\trois{simple-type}
Any~$x\in D(F)$ defines a type on~$D$, and one can identify that way~$D(F)$ with a subset of~$S(D)$, whose elements
will be called {\em simple}
types.

\trois{type-functoriality}
{\em Functoriality}.
Let $f$ be an $F$-definable natural transformation from $D$ to another $F$-definable
functor $D'$ and let $x\in S(D)$. If $t$ is a representative of $x$ then the $F$-equivalence class
of $f(t)$ only depends on $x$, and not on $t$; one thus gets a well-defined type on $D'$, 
that will usually 
be denoted by $f(x)$. 

\medskip
If~$\Delta$ is an~$F$-definable sub-functor of~$D$, the induced
map $S(\Delta) \to S(D)$ identifies
$S(\Delta)$ with the subset of~$S(D)$
consisting of types admitting a representative~$(L,x)$ with~$x\in \Delta(L)$ 
(it will then be the case for {\em all}
its representatives,
by the very definition of~$F$-equivalence); we simply say that such a type~{\em lies on~$\Delta$.}

\trois{type-extension}
Let~$L\in \mathsf M_F$, and let~$D_L$ be the restriction of~$D$
to~$\mathsf M_L$. Let~$t$ be a type on~$D_L$.  
Any representative of~$t$ defines 
a type on~$D$, which only depends on~$t$
(indeed, $L$-equivalence is stronger than~$F$-equivalence); one thus
has a natural 
{\em restriction}
map~$S(D_L)\to S(D)$.

\trois{ultra-filters}
{\em Remark.} Let~$x\in S(D)$. Let~$\mathscr U_x$
be the set of subsets of~$D(F)$ which are of the kind~$\Delta(F)$, where~$\Delta$
is an~$F$ definable sub-functor of~$D$ on which~$x$ lies. One can rephrase the celebrated
compactness theorem of model theory by saying that~$x\mapsto \mathscr U_x$ establishes a bijection
between~$S(D)$ and the set of ultra-filters of~$F$-definable subsets of~$D(F)$.

\trois{def-types}
{\bf Definable types}. 
Let~$t\in S(D)$. By definition, the type $t$ is determined by the 
data of all~$F$-definable
sub-functors of~$D$ on which it lies. We will say that~$t$ is {\em $F$-definable} if, roughly speaking, the following holds: 
{\em for every~$F$-definable family of~$F$-definable sub-functors of~$D$, the set of parameters for which~$t$ lies on
the corresponding~$F$-definable sub-functor is itself~$F$-definable.}

\medskip
Let us be more precise. If~$D'$ is an~$F$-definable functor, and if~$\Delta$ is
an~$F$-definable sub-functor
of~$D\times D'$, then for every~$L\in \mathsf M_F$ 
and every~$x\in D'(L)$, the fiber~$\Delta_x$ of~$\Delta_L$ over~$x$ is an~$L$-definable
sub-functor~$\Delta_x$ of~$D_L$. 

We will say that~$t$ is~$F$-definable if for every such~$(D',\Delta)$, the set of~$x\in D'(F)$ such that~$t$ lies
on~$\Delta_x$ is a definable subset of~$D'(F)$.

\trois{canon-extension}
{\bf Canonical extension of an $F$-definable type.} 
If~$t$ is an~$F$-definable type on~$D$, then it admits for every~$L\in \mathsf M_F$ a canonical 
$L$-definable pre-image~$t_L$ on~$S(D_L)$, which is called
the {\em canonical extension}
of~$t$.

\medskip
Roughly speaking,~$t_L$ is defined by the same formulas as~$t$. 
This means the following. Let~$\Sigma$ be an~$L$-definable sub-functor of~$D_L$. Considering the coefficients of the formulas
that define~$\Sigma$ as parameters, we see that there exist an~$F$-definable functor~$D'$, an~$F$-definable
sub-functor~$\Delta$ of~$D\times D'$, and a point~$y\in D'(L)$ such that~$\Sigma=\Delta_y$. 

Now since~$t$ is~$F$-definable, there exists an~$F$-definable sub-functor~$E$ of~$D'$ such that for every~$x\in D'(F)$, 
the type~$t$ lies on~$\Delta_x$ if and only if~$x\in E(F)$. Then the type~$t_L$ belongs to~$\Sigma=\Delta_y$ if and only if~$y\in E(L)$.

\trois{ortho-gamma}
{\bf Orthogonality to~$\Gamma$.} 
Let~$t$ be a type on~$D$. We will say that~$t$ is {\em orthogonal to~$\Gamma$}
if it is~$F$-definable and if
for every~$F$-definable natural transformation~$f: D\to \Gamma_{0,F}$, the image~$f(t)$, which is {\em a priori}
a type on~$\Gamma_{0,F}$, is a {\em a simple} type, that is, belongs to~$\Gamma_0(F)=|F|$. 
If this is the case, $t_L$ remains orthogonal to~$\Gamma$ for every~$L\in \mathsf M_F$.

\medskip
It follows from the definitions that any simple type on~$D$ is orthogonal to~$\Gamma$. 

\trois{functoriality-typedef} {\bf Functoriality of those definition}. If $f\colon D\to D'$ is 
an $F$-definable natural transformations, then for every
type $t\in S(D)$ the image $f(t)$ is $F$-definable (resp. orthogonal to $\Gamma$)
as soon as $t$ is $F$-definable (resp. orthogonal to $\Gamma$). 

\deux{cas-alg-var}
{\bf A very important example : the case of an algebraic variety.} 
Let~$Y$ be an algebraic variety over~$F$. The purpose of this
paragraph is to get an explicit description of~$S(Y)$; as we will see,
this descriptions looks like that of the Berkovich analytification of an algebraic 
variety (\ref{analytify}). 

\trois{alg-var-gen}
{\bf The general case}. Let~$t$ be a type on~$Y$. Let~$(L,x)$
be a representative of~$t$. The point~$x$ of
$Y(L)$ 
induces a {\em scheme-theoretic} point~$y$
of~$Y$ and a valuation~$|.|_x$
on the residue field~$F(y)$, extending that of~$F$; the
data of the point~$y$ and of (the equivalence class of)
the valuation~$|.|_x$ only depend on~$t$, 
and not on the choice of~$(L,x)$. 

Conversely, if~$y$ is a schematic
point of~$Y$, any valuation on~$F(y)$ extending that of~$F$ is induced by an isometric embedding~$F(y)\hookrightarrow L$
for some~$L\in \mathsf M_F$, hence arises from some type~$t$ on~$Y$. One gets that way a bijection
between~$S(Y)$ and the {\em valuative spectrum}
of~$Y$, that is, the set of couples~$(y,|.|)$ with~$y$ a scheme-theoretic point on~$Y$ and $|.|$
a valuation on~$F(y)$ extending that of~$F$ (considered up to equivalence). Note that there is a natural 
map from $S(Y)$ to $Y$, sending a couple $(y,|.|)$ to $y$. 

\trois{aff-case}
{\bf The affine case.}
Let us assume now that~$Y=\spec B$, let~$y\in Y$ and let~$|.|$ be a valuation on~$F(y)$ extending
that of~$F$. By composition with the evaluation map~$f\mapsto f(y)$, we get a valuation~$\phi : B\to |F(y)|$ (the definition of a valuation on a ring is 
{\em mutatis mutandis}
the same as on a field -- but be aware that it may have a non-trivial kernel); this valuation extends that of~$F$. 

Conversely, let~$\phi$ be a valuation on~$B$ extending that of~$F$. 
Its kernel corresponds to a point~$y$
of~$Y$, and~$\phi$ is the composition of the evaluation map at~$y$ 
and of a valuation~$|.|$ on~$F(y)$ extending that of~$F$. 

These constructions provide a bijection between
the valuative spectrum of~$Y$ and the set of (equivalence classes of) valuations on~$B$ extending that of~$F$; 
eventually, we get a bijection between~$S(Y)$ and the set of valuations on~$B$ extending that of~$F$. Modulo this
bijection, the natural map $S(Y)\to Y$ sends a valuation to its kernel.

\trois{interpretation-sy}
{\bf Interpretation of some properties}. Let $t\in S(Y)$
let $\phi$ be the corresponding valuation on~$B$, and let $\got p$
be the kernel of $\phi$ (that is, its image on $Y=\spec B$). 
Let $\got G$ be the subset of $B^2$ consisting of couples 
$(b,b')$ such that $\phi(b)\leq \phi(b')$.  

\medskip
{\em Interpretation of definability.}
The type~$t$ is~$F$-definable if and only if for every finite~$F$-dimensional
subspace~$E$ of~$B$ the following holds: 

- the intersection $E\cap \mathfrak p$ is an $F$-definable subset of $E$ ; 

- the intersection $(E\times E)\cap \got G$ is an $F$-definable subset of $E\times E$.

\medskip
{\em Interpretation of orthogonality to~$\Gamma$.}
The type~$t$ is orthogonal to~$\Gamma$ if and only it is~$F$-definable, and if~$\phi$ takes its values in~$|F|$. This
is equivalent to require that~$\phi$ takes its values in~$|F|$ and that~$\phi_{|E}: E\to |F|$ is~$F$-definable for every
finite~$F$-dimensional subspace~$E$ of~$B$. 

\medskip
{\em Interpretation of the canonical extension.}
Assume that
the type $t$ is~$F$-definable and let~$L\in \mathsf M_F$; we want to describe
the valuation~$\phi_L$
on~$B\otimes_F L$ that corresponds to the canonical extension~$t_L$
(\ref{canon-extension}). 
It is equivalent to describe the kernel $\got p_L$ of $\phi_L$ and 
the subset $\got G_L$ of $(B\otimes_FL)^2$ consisting of couples
$(b,b')$ such that $\phi_L(b)\leq \phi_L(b')$. 
It is sufficient to describe the intersection $\got p_L$ (resp. $\got G_L$)
with~$E\otimes_k L$ (resp. with $(E\otimes_k L)^2$)
for any finite dimensional $F$-vector subspace~$E$ of~$B$. 

Let us fix such an~$E$. Since~$t$ is definable, there exist a unique $k$-definable sub-functor~$D$ of~$\underline E:=\Lambda \mapsto E\otimes_F \Lambda$
such that $E\cap \got p=D(F)$
and a unique $k$-definable sub-functor $D'$ of $\underline E^2$ such that
$(E\times E)\cap \got G=D'(F)$. One then has
$$\got p_L\cap (E\otimes_FL)=D(L)\;\;\;{\rm and}\; \;\;\got G_L \cap (E\otimes_F L)^2=D'(L).$$

\medskip
Let us assume moreover that~$t$ is orthogonal to~$\Gamma$. In this situation,
$\phi$ induces
a~$k$-definable map~$E\to |F|$. This definable map comes from a unique~$k$-definable transformation~$\underline E \to \Gamma_0$. 
By evaluating it on~$L$, one gets an~$L$-definable map~$E\otimes_F L \to |L|$, which coincides with the restriction of~$\phi_L$. 

\deux{def-vhat}
{\bf The fundamental definition (\cite{hl} \S 3).}
Let~$V$ be a~$k$-definable functor. The {\em stable completion}\footnote{This is named after the model-theoretic notion of
{\em stability} 
which plays a key role in Hrushovski and Loers's work
 ({\em cf.} \cite{hl} \S 2.4, 2.6 and 2.9), but is far beyond our scope and
 will thus in some sense 
remain
hidden in this text.} of $V$
is the functor $$\hat V : \mathsf M \to \mathsf{Sets}$$
defined as follows: 

\medskip
$\bullet$ if~$F\in \mathsf M$, then~$\hat V(F)$ is the set of types on~$V_F$ that are orthogonal to~$\Gamma$; 

$\bullet$ if~$L\in \mathsf M_F$, 
the arrow~$\hat V(F)\to \hat V(L)$ is the embedding that sends a type~$t$ to its canonical extension~$t_L$
(\ref{canon-extension}). 

\deux{property-vhat}
{\bf Basic properties and first examples}.

\trois{funct-vhat}
The formation of~$\hat V$ is functorial in~$V$ with respect to~$k$-definable maps (\ref{functoriality-typedef}). 

\medskip
\trois{simple-points}
Since any simple type
(\ref{simple-type})
is orthogonal to~$\Gamma$, one has a natural embedding of functors~$V\hookrightarrow \hat V$. 
We will therefore identify~$V$ with a sub-functor of~$\hat V$. For every~$F\in \mathsf M$, the points
of~$\hat V$
that belong to~$V(F)$ will be called
{\em simple}
points. 

\trois{polyhedron-hat} {\bf The stable completion of a polyhedron}.  
Let~$F\in \mathsf M$. It follows immediately from the definition 
that a type on~$\Gamma_{0,F}$ is orthogonal to~$\Gamma$ if and only if it is simple. 
In other words, $\hat {\Gamma_0}=\Gamma_0$. 
This fact extends to $k$-definable sub-functors of~$\Gamma_0^n$ (\ref{definable-pl}): if~$V$ is such a functor
then the natural embedding~$V\hookrightarrow \hat V$
is a bijection. 

\deux{variety-hat}
{\bf The stable completion of a definable sub-functor of an algebraic variety}. 
Let~$X$ be a $k$-scheme of finite type, and let~$V$ be a~$k$-definable sub-functor of~$X$.

\trois{x=v=spec} {\em Assume that $V=X=\spec A$}. 
From~\ref{aff-case}
and~\ref{interpretation-sy}
we get the following description of~$\hat X$. Let~$F\in \mathsf M$. 
The set $\hat X(F)$ is the set of 
valuations $$\phi : A\otimes_k F\to |F|$$
extending that of~$F$ and 
such that for every finite dimensional~$F$-vector space~$E$ of~$A\otimes_k F$, the restriction $\phi_{|E} : E\to |F|$ is~$F$-definable.
Let~$L\in \mathsf M_F$ and let~$\phi\in \hat X(F)$. For every finite dimensional~$F$-vector subspace~$E$ of~$A\otimes_k F$, the~$F$-definable map
$\phi_{|E} : E\to |F|$ arises from a unique~$F$-definable natural transformation~$\Phi_E : \underline E \to \Gamma_0$, which itself gives rise to a~$L$-definable map
~$\Phi_E(L) : E\otimes_F L\to |L|$. By gluing the maps $\Phi_E(L)$'s for $E$ going through the set of all finite $F$-dimensional vector subspaces of $A$
we get a valuation ~$A\otimes_k L \to |L|$ extending that of~$L$; this is precisely 
the image~$\phi_L$ of~$\phi$ under the natural embedding~$\hat X(F)\hookrightarrow \hat X(L)$. 
Roughly speaking, $\phi_L$ is ``defined by the same formulas as~$\phi$.''

We thus see that Hrushovski and Loeser
mimic in some sense Berkovich's construction
(\ref{analytify}, but with 
a model-theoretic and
definable flavour. 

\trois{speca-vqque}
{\em Assume that $X=\spec A$, and that~ $V$ is defined by a boolean combination
of inequalities of the form $|f|\Join \lambda |g|$ (with~$f$ and~$g$ in~$A$ and~$\lambda\in |k|$)}. 
The functor~$\hat V$ is then the sub-functor of~$\hat X$ consisting of the semi-norms~$\phi$
satisfying the same combination of inequalities. 

\trois{vhat-general} In general, one defines ~$\hat V$
by performing the above constructions locally and gluing them.

\deux{vhat-topology} {\bf A topology on~$\hat V(F)$.}
Let~$V$ be a~$k$-definable
functor and let~$F\in \mathsf M$. We are going to define a topology on~$\hat V(F)$
in some particular cases; all those constructions will we based upon the {\em order topology}
on~$|F|$.

\medskip
\trois{topology-polyhedron}
If~$V$ is a~$k$-definable sub-functor of~$\Gamma_0^n$
then~$\hat V(F)=V(F)$ is endowed
with the topology induced from the product topology on~$\Gamma_0(F)=|F|^n$. 
In particular if $I$ is a~ $k$-definable interval it inherits a topology. 

\trois{topology-generalized}
If~$V$ is a~$k$-definable generalized segment, then $\hat V(F)=V(F)$
inherits a topology using the above construction and the quotient topology.

\trois{topology-variety}
If~$V$ is a~$k$-scheme of finite type, then $\hat V(F)$ is given the coarsest topology such that: 

\medskip
$\bullet$ for every Zariski-open subset $U$ of $V$, the subset $\hat U(F)$ of $\hat V(F)$ is open; 

$\bullet$ for every affine open subset $U=\spec A$ of $V$ and every $a\in A\otimes_k F$, the map
$\hat U(F)\to |F|$ obtained by applying valuations to $a$ is continuous.

\trois{topology-def-variety}
If $V$ is a $k$-definable subfunctor of a $k$-scheme of finite type $X$ then $\hat V(F)$
is endowed with the topology induced from that of $\hat X(F)$, defined at~\ref{topology-variety} above.

\medskip
\deux{topology-comments}
{\bf Comments.} 

\trois{presentation}
Let~$V$ be a $k$-definable functor which is a sub-functor of $\Gamma_0^n$, a $k$-definable generalized segment or a 
sub-functor of a $k$-scheme of finite type.  For a
given $F\in \mathsf M$ the set $\hat V(F)$
inherits a topology as described above at~\ref{vhat-topology} {\em et sq.} 

We emphasize that this topology does not depend only on the abstract $k$-definable functor $V$, but also 
on its given presentation. 
For example, let $i \colon (\spec k)\coprod \mathbb G_{{\rm m},k}\to \mathbb A^1_k$
be the morphism induced by the closed immersion of the origin and the open immersion of $ \mathbb G_{{\rm m},k}$
into $\mathbb A^1_k$. The natural transformation induced by $i$ is then a $k$-definable isomorphism, because for every 
$F\in \mathsf M$ the induced map $\{0\}\coprod F^\times \to F$ is bijective. But it is not a homeorphism: $\{0\}$ is open 
in the left-hand side, and not in the right-hand side. 

\trois{point-topology}
Let~$V$ be a $k$-definable sub-functor of some
$k$-scheme of finite type $X$. The valuation on $F$ defines a topology on $X(F)$; if $V(F)$ is given the induced topology the inclusion
$V(F)\subset \hat V(F)$ is then a {\em topological embedding with dense image}. 

%

\trois{topology-ground}
Let~$F\in \mathsf M$ and let~$L\in \mathsf M_F$. Be aware that in general, the topology on~$|F|$ 
induced by the order topology on~$|L|$ is
{\em not} the order topology on~$|F|$; it is finer. 
For instance, 
assume that there exists~$\omega \in |L|$ such that~$1<\omega$ and
such that there is no element~$x\in |F|$ with~$1<x\leq \omega$ (in other words,
$\omega$ is upper infinitely close to 1
with respect to~$|F|$). Then for every~$x\in |F^*|$ the singleton~$\{x\}$ is equal to 
$$\{y\in |F|, \;\;\omega^{-1} x < y< \omega x\}.$$ It is therefore
open for the topology induced by the order topology on~$|L|$; hence the latter
induces the discrete topology on~$|F^*|$. 

Since all the topologies we have considered above ultimately rely on the order topology, 
this phenomenon also holds for them. That is, let~$V$
be as in~\ref{presentation}. The natural
embedding~$\hat V(F)\hookrightarrow \hat V(L)$ is {\em not}
continuous in general; the topology on~$\hat V(F)$ is coarser than the topology induced
from that of~$\hat V(L)$.

\trois{continuous-maps}
Let~$V$ and~$W$ be two $k$-definable functors, each of which 
being of one of the 
form considered in~\ref{presentation}. We will say that a $k$-definable
natural transformation~$f: \hat V\to \hat W$
is {\em continuous} if~$f(F): \hat  V(F)\to \hat W(F)$ is continuous for every~$F\in \mathsf M$; we then 
define
in an analogous way the fact for~$f$ to be a homeomorphism, or to induce a homeomorphism
between~$\hat V$ and a sub-functor of~$\hat W$, etc.

\deux{example-1}{\bf The affine and projective lines (\cite{hl}, example 3.2.1).} 
Let~$F\in \mathsf M$.

\trois{def-etar}
For every~$a\in F$ and~$r\in |F|$, the map~$\eta_{a,r,F}$
from~$F[T]$ to~$|F|$ that sends
$\sum a_i (T-a)^i$
to~$\max |a_i|\cdot r^i$
is a valuation on~$F[T]$ whose associated type
belongs to~$\widehat{\Aff^1_k}(F)$. 
Note that for every~$a\in F$, the semi-norm~$\eta_{a,0,F}$ is nothing but~$P\mapsto |P(a)|$; hence
it corresponds to the simple point~$a\in F=\Aff^1_k(F)$. 
If~$L\in \mathsf M_L$, the natural embedding
$\widehat{\Aff^1_k}(F)\hookrightarrow \widehat{\Aff^1_k}(L)$ 
is induced by the map that 
sends the semi-norm~$\eta_{a,r,F}$ (for given~$a\in F$ and~$r\in |F|$) to the semi-norm
on~$L[T]$ ``that is defined by the same formulas'', that is, to~$\eta_{a,r,L}$. 

\medskip
\trois{rem-etar}
{\em Remark.}
If the ground field~$F$ is clear from the context, we will 
sometimes write simply~$\eta_{a,r}$ instead of~$\eta_{a,r,F}$. 

\trois{balls-a1}
For every~$(a,b)\in F^2$ and~$(r,s)\in |F|^2$, 
an easy computation shows that 
the valuations $\eta_{a,r}$ and~$\eta_{b,s}$ are equal 
if and only if~$r=s$ and~$|a-b|\leq r$, that is, if
and only if the closed balls~$B(a,r)$ and~$B(b,s)$ of~$F$ are equal. 
Moreover, one can prove that every valuation
belonging to~$\hat{\Aff^1_k}(F)$ is of the form~$\eta_{a,r}$ 
for suitable~$(a,r)\in F\times |F|$. Therefore we get a 
functorial
{\em bijection} between $\widehat{\Aff^1_k}(F)$
and the set of closed balls of~$F$; {\em it then follows from~\ref{ball-definable}
that the functor
$\widehat{\Aff^1_k}$ is~$k$-definable}.

\trois{p1-def}
The functor~$\widehat{\PP^1_k }$ is simply obtained by adjoining 
the simple point~$\infty$ to~$\widehat{\Aff^1_k}$; hence it is $k$-definable too. 

%
%
%

\deux{pro-def}{\bf Definitions.}
In order to be able to state fundamental results about the stable completions, we need some 
definitions.

\trois{dimensions-definable}
Let $X$ be a $k$-scheme of finite type and let 
$V$ be a $k$-definable subfunctor of $X$. Let $F\in \mathsf M$. The set
of integers $n$ such that there exists an $F$-definable
injection $(F\zero)^n\hookrightarrow V(F)$ does not depend on $F$
(this follows from quantifier elimination)
and is bounded by $\dim X$. Its supremum is called the {\em dimension of $V$}; 
this is a non-negative integer $\leq \dim X$ if~$V\neq \varnothing$, and $(-\infty)$ otherwise.

\medskip
Let us mention for further use that there is also a notion of dimension for $k$-definable sub-functors of $\Gamma_0^n$. 
Let $D$ be such a sub-functor and let $F\in \mathsf M$. The set of integers $m$ such that there exists 
elements $a_1,b_1,\ldots, a_m, b_m$ in $|F^\times|$ with $a_i<b_i$ for every $i$ and a $k$-definable injection 
$$[a_1;b_1]\times\ldots \times [a_m;b_m]\hookrightarrow D$$ does not depend on $F$. Its supremum is called the dimension of $D$. 
This is a non-negative integer $\leq n$ if~$D\neq \varnothing$, and $(-\infty)$ otherwise. 

\trois{def-prodefinable}
A functor from $\mathsf M$ to $\mathsf{Sets}$ is said to be {\em pro-$k$-definable} if it is isomorphic to a projective limit of~$k$-definable functors. 
This being said, let us quickly mention two issues we will not discuss in full detail
(the interested reader may refer to \cite{hl}, \S 2.2). 

1) If one wants this isomorphism
with a projective limit to be canonical, one has to be implicitly 
given a definition of this functor ``in families''
(compare to~\ref{abstract-embedable}); this is the case here: there is a natural notion of
a $k$-definable family of types orthogonal to~$\Gamma$. 

2) The set of indices of the projective system has to be of small enough (infinite) cardinality; in the cases we consider in this
paper it can be taken
to be countable.

\trois{prodef-map}
Let $V=\projlim V_i$ and $W=\projlim W_j$ be two pro-$k$-definable functors. For every $(i,j)$, let
$D_{ij}$ be the set of all $k$-definable natural transformations from $V_i$ to $W_j$. 
A natural transformation from $V$ to $W$ is said to be {\em pro-$k$-definable}
if it belongs to 
$$\lim\limits_{\stackrel\longleftarrow j}\;
\left(\lim\limits_{\stackrel \longrightarrow i}\; D_{ij}\right).$$

\trois{strict-prodef}
A functor $V$ from $\mathsf M$ to $\mathsf {Sets}$ 
is said to be {\em strictly}
pro-$k$-definable if
for every pro-$k$-definable natural transformation~$f$
from~$V$ to a~$k$-definable functor~$D$, the direct image functor~$f(V)$ is
a~$k$-definable sub-functor of~$D$ (it is {\em a priori}
only a (possibly infinite) intersection of~$k$-definable sub-functors). 

\trois{rel-prodef}
Let $V$ be a pro-$k$-definable functor from
$\mathsf M$ to $\mathsf{Sets}$. A sub-functor $D$ of $V$
is said to be {\em relatively $k$-definable} if there exists
a $k$-definable functor $\Delta$, 
a pro-$k$-definable natural transformation~$f: V\to \Delta$
and a $k$-definable sub-functor~$\Delta'$ of~$\Delta$
such that~$D=f^{-1}(\Delta')$. 

\medskip
Let $F\in \mathsf M$. 
By evaluating those functors at $F$-points, one gets the notion of a {\em relatively~$k$-definable}
subset of~$V(F)$, and more generally of a~{\em relatively~$F_0$-definable}
subset of~$V(F)$ for~$F_0$ any valued field lying between~$k$ and~$F$
(compare to~\ref{scalar-change}). 

\medskip
Replacing~$F$-definable sub-functors by
relatively~$F$-definable ones, 
in the definitions of~\ref{type}
{\em et sq.}
one gets the notion of a type
on~$V_F$, and of an~$F$-definable type
on~$V_F$.
There is a natural bijection between the set of types on~$V_F$
and that of ultra-filters of relatively~$F$-definable subsets of~$V(F)$
(compare to remark~\ref{ultra-filters}).

\deux{hatv-prodef}
{\bf Theorem (\cite{hl}, lemma 2.5.1, th. 3.1.1, th. 7.1.1 and rem. 7.1.3)}. {\em Let $V$ be a $k$-definable sub-functor of a $k$-scheme 
of finite type. The stable completion~$\hat V$ is strictly 
pro-$k$-definable, 
and is $k$-definable if and only if~$\dim V\leq 1$}. 

\trois{quick-words-prodefinability}
The pro-$k$-definability of~$\hat V$ comes from general arguments of model-theory, which 
hold in a very general context, that is, not only in the theory of valued fields. 
Its {\em strict}
pro-$k$-definability 
comes from model-theoretic properties of the theory of algebraically closed fields, 
used at the ``residue field'' level. 
 
 \trois{def-vhat-dim1}
 The
 definability of~$\hat V$ in dimension~$1$
is far more specific. It ultimately relies
on the 
Riemann-Roch theorem for curves,
through 
the following consequence of the latter: if~$X$ is a projective, irreducible, smooth curve of genus~$g$ over
an algebraically closed field~$F$, the group~$F(X)^\times$ is generated by rational functions with
at most~$g+1$ poles (counted with multiplicities).

\deux{compactness}
{\bf Definable compactness.} Let $V$ be a $k$-definable sub-functor
of $\Gamma_0^n$, a $k$-definable generalized segment or a $k$-definable
subscheme of a scheme of finite type. The functor $\hat V$ is then pro-$k$-definable: indeed, 
in the first two cases it is equal to $V$, hence even $k$-definable; in the third case, this is th.~\ref{hatv-prodef}
above. 

\trois{limit-type} 
We have defined at~\ref{vhat-topology}
{\em et sq.} for
every~$F\in \mathsf M$ a topology
on~$\hat V(F)$. 
By construction, there exists for every~$F\in \mathsf M$
a set~$\mathbb O_F$ of relatively~$F$-definable sub-functors of~$\widehat {V_F}$  such that:

\medskip
$\bullet$ the set~$\{D(F)\}_{D\in \mathbb O_F}$
is a basis of open subsets of~$\hat V(F)$; 

$\bullet$ for every~$D\in \mathbb O_F$ and every~$L\in \maths M_F$, the sub-functor~$D_L$
of~$\widehat {V_L}$ belongs to~$\mathbb O_L$ (in particular, $D(L)$ is an open subset
of~$\hat V(L)$). 

\medskip
Let~$F\in \mathsf M$ and let~$t$ be a type on~$\widehat{V_F}$
(this makes sense in view of~\ref{rel-prodef}
since $\hat V$ is pro-$k$-definable by th.~\ref{hatv-prodef}). 
We say that a point~$x\in \hat V(F)$ is a {\em limit}
of~$t$ if the following holds: {\em for every~$D\in \mathbb O_F$ such that~$x\in D(F)$,
the type~$t$ lies on~$D$.} If~$X$ is separated~$\hat V(F)$ is Hausdorff, 
and the limit of a type is then unique provided it exists.  

\trois{def-compactness}
We will say that~$\hat V$ is
{\em definably compact}
if
for every~$F\in \mathsf M$, every
{\em $F$-definable}
type on~$\widehat {V_F}$ has a unique
limit in~$\hat V(F)$. 

\trois{rem-ultrafilters-compact}
{\em Remarks.}
Viewing types as ultra-filters, 
we may rephrase the above definition by saying that $\hat V$
is definably compact if for every~$F\in \mathsf M$, every~$F$-definable 
ultra-filter of relatively~$F$-definable subsets of~$\hat V(F)$
has a unique limit on $\hat V(F)$.  

\deux{avatar-compactness}
The following propositions provides evidence for 
definable compactness being the right model-theoretic analogue
of usual compactness; they are particular caseds
of some general results
proved (or at least stated with references) in \cite{hl}, \S 4.1 and 4.2. 

\trois{polyhedron-compact}
{\bf Proposition.}
{\em Let $V$ be a $k$-definable sub-functor of $\Gamma_0^n$. 
It is definably compact if and only if it is contained in $[0;R]^n$ for some $R\in |k|^{\QQ}$
and can be defined by conjunction
and disjunction of {\em non-strict}
monomial inequalities}.  

\trois{rem-gensem-compact}
{\em Remark}. Any $k$-definable generalized segment
is definably compact: this can be proved either directly, or by writing
down an explicit $k$-definable {\em homeomorphism}
between such a segment and a functor $V$ as in prop.~\ref{polyhedron-compact}
above. 

\trois{projective-compact}
{\bf Proposition}. {\em Let $X$ be a $k$-scheme of finite type. 

\medskip
1) The stable completion $\hat X$ is definably compact if and only if $X$ is proper. 

2) Assume that $X$ is proper, let $(X_i)$ be an affine covering of $X$
and let $V$ be a $k$-definable sub-functor of $X$ such that
$V\cap X_i$ can be defined for every~$i$
by conjunction and disjunction of {\em non-strict} inequalities. The stable
completion~$\hat V$ is then definably compact.}

\section{Homotopy type of~$\hat V$ and links with Berkovich spaces}

\deux{link-berko}
{\bf The link between stable completions and Berkovich spaces}. 
For this paragraph, we assume that
the valuation~$|.|$ of~$k$
takes real values, and that~$k$ is complete. 

\trois{semi-algebraic}
Let~$X$ be an algebraic variety over~$k$, and let~$V$ be a~$k$-definable sub-functor of~$X$. 
The inequalities that define~$V$ also define a strict
semi-algebraic subset~$V^{\rm an}$ of $X\an$ (for the definition of such a subset
see~\ref{semi-alg});  as the notation suggests, $V\an$
only
depends on the sub-functor $V$ of $X$, and not on its chosen description
(this follows from quantifier elimination). 
Any strict semi-algebraic subset of~$X^{\rm an}$ is of that kind.

\trois{def-pi}
Let~$F\in \mathsf M$ be such that~$|F|\subset \RR_+$. 
Any point of~$\hat V(F)$
can be interpreted in a suitable affine chart as a valuation
with values in~$|F|\subset \RR_+$; it thus induces
a point 
of~$V^{\rm an}$. 
We get that way a map~$\pi : \hat V(F)\to V^{\rm an}$
which is
continuous by the very definitions of the topologies
involved. 
One can say a bit more about~$\pi$ in some particular cases; the following
facts are proved in \cite{hl}, \S 14.

\medskip
\trois{case-k-F}
{\em We assume that $k$ is algebraically closed, non-trivially valued and that $F=k$.}
The map~$\pi$ then induces a homeomorphism
onto its image.  

\medskip
Let us describe this image for~$V=X=\Aff^1_k$. By the explicit
description of~$\widehat {\Aff^1_k}$, it consists precisely of the 
set of points~$\eta_{a,r}$ with~$a\in k$ and~$r\in |k|$, that is, of the set
of points of type 1 and 2 according to Berkovich's classification (see \cite{brk1}, chapter 2 and 4). 

This generalizes as follows: if~$X$ is any curve, then~$\pi(\hat V(k))$ is precisely the set of
points of~$V^{\rm an}$ that are of type 1 or 2.

\medskip
{\em Comments.} The fact that Hrushovski and Loeser's theory, 
which focuses on definability, only sees points of type 1 and type 2, encodes the following
(quite vague) phenomenon: {\em when one deals only with algebraic curves and scalars belonging
to the value groups of the ground field, points of type 1 and type 2 are the only ones at which 
something interesting may happen.} 

Of course, considering all Berkovich points is useful 
because it ensures 
good topological properties (like pathwise connectedness and compactness if one starts from 
a projective, connected variety). Those properties are lost
when one only considers points of type 1 and type 2;
Hrushovski and Loeser
remedy it
by introducing the corresponding model-theoretic properties,
like definable compactness, or definable 
arcwise 
connectedness
(see below for the latter). 

\trois{case-f-maximally}
{\em We assume that $|F|=\RR_+$ and that
$F$ is
{\em maximally complete} ({\em i.e.}
it does not admit
any proper valued extension with the same value group
and the same residue field)}. 
In this case (without any particular assumption on $k$)
the continuous map~$\pi$ is a proper surjection. 

\trois{k-F-maxi} {\em We assume that $|k|=\RR_+$, that $k$
is algebraically closed and maximally complete, and that $F=k$.} Fitting together
\ref{case-k-F} and~\ref{case-f-maximally}
we see that~$\pi$
establishes then a homeomorphism~$\hat V(k)\simeq V^{\rm an}$. 

\medskip
Suppose for the sake of simplicity that~$V=X=\spec A$. By definition, 
a point of~$X^{\rm an}$ is a valuation~$\phi : A\to \RR_+=|k|$ extending that of~$k$. 
Therefore the latter assertion says that for any
finite dimensional~$k$-vector space~$E$ of~$A$ the restriction~$\phi_{|E}$
is 
{\em automatically}~$k$-definable, because of the
maximal completeness of~$k$. This ``automatic definability'' result comes
from the previous work~\cite{hhm}
by Haskell, Hrushovski and Macpherson about
the model theory of valued fields (which is intensively used by Hrushovski
and Loeser throughout their paper); the reader who would like to have 
a more effective version of it with a somehow explicit
description 
of the formulas describing~$\phi_{|E}$
may refer to the recent work~\cite{poineautypes}
by Poineau.

\deux{presentation-hl-work}
Hrushovski and Loeser therefore work
{\em purely inside
the ``hat'' world}; they transfer thereafter their results to the Berkovich setting (in \cite{hl}, \S 14)
using the
aforementioned map~$\pi$ and its good properties in the maximally complete case
(\ref{case-f-maximally}, \ref{k-F-maxi}). 

\medskip
For that reason, {\em we will not deal with Berkovich spaces
anymore}. We are now going to state the ``hat-world'' avatar of~\ref{def-ret-hl}
(this is th.~\ref{main-theo} below)
and sketch its proof
very roughly. 

\medskip
In fact, th.~\ref{main-theo}
is more precisely the ``hat world" avatar
of a {\em weakened} version of~\ref{def-ret-hl}, consisting of the 
same statement
for {\em strict} (instead of general)
semi-algebraic subsets of Berkovich spaces. 
But this is not
too serious a restriction: indeed, Hrushovski and Loeser prove their
avatar of~\ref{def-ret-hl}
by first reducing it to the strict case
through a nice trick,
consisting more or less in ``seeing 
bad scalars as
extra
parameters'' (\cite{hl}, beginning of \S 11.2; the reader may also 
refer to
\cite{bourb}, \S 4.3 for
detailed explanations).

\deux{main-theo} {\bf Theorem (simplified version of th.
11.1.1 of
\cite{hl}).}
{\em Let~$X$ be a {\em quasi-projective}~$k$-variety, 
and let~$V$ be a~$k$-definable sub-functor of~$X$. Let~$G$
be a finite group acting on~$X$
and stabilizing~$V$, 
and let~$E$
be a finite set of~$k$-definable transformations
from~$V$
to~$\Gamma_0$ (if~$f\in E$, 
we still denote by~$f$ the induced natural transformation
$\hat V\to \widehat{\Gamma_0}=\Gamma_0$; if~$g\in G$, we still 
denote by~$g$ the induced 
automorphism of~$\hat V$).There exists: 

\medskip
\begin{itemize} 

\item[$\bullet$] a~$k$-definable generalized segment~$I$, with endpoints~$o$ and~$e$; 

\item[$\bullet$] a~$k$-definable sub-functor~$S$ of~$\hat V$ (where $S$ stands for {\em skeleton}); 

\item[$\bullet$] a $k$-definable sub-functor $P$ of $\Gamma_0^m$ (for some $m$) of dimension~$\leq \dim X$; 

\item[$\bullet$] a~$k'$-definable homeomorphism~$S\simeq P$, for a suitable finite extension~$k'$ of~$k$
inside~$\overline k$; 

\item[$\bullet$] a continous~$k$-definable map~$h: I\times \hat V \to S$
satisfying the following properties
for every~$F\in \mathsf M$, every~$x\in \hat V(F)$, every~$t\in I(F)$, every~$g\in G$
and every~$f\in E$: 

\begin{itemize}
\item [$\diamond$] $h(o,x)=x$, and~$h(e,x)\in S(F)$; 

\item[$\diamond$] $h(t,x)=x$ if~$x\in S(F)$; 

\item[$\diamond$] $h(e,h(t,x))=h(e,x)$;

\item[$\diamond$] $f(h(t,x))=f(x)$; 

\item[$\diamond$] $g(h(t,x))=h(t,g(x))$. 

\end{itemize}
\end{itemize} 

 }

\deux{comments-main-theo}
{\em Comments.}

\trois{twist-poly}
We will refer to the existence of~$k'$, of~$P$,
and of the~$k'$-definable homeomorphism~$S\simeq P$
by saying that~$S$ is a {\em tropical}.
The finite extension~$k'$ of~$k$
cannot be avoided; indeed, it
reflects the fact that the Galois action on the homotopy type
of~$\hat V$ is non necessarily trivial. In the Berkovich language, think
of the~$\QQ_3$-elliptic curve~$E: y^2=x(x-1)(x-3)$. The analytic curve~$E_{\QQ_3(i)}^{\rm an}$ 
admits a Galois-equivariant deformation retraction to a circle, on which
the conjugation exchanges two
half-circles; it descends to a deformation retraction of~$E^{\rm an}$
to a compact interval. 

\trois{def-homotopy}
We will sum up
the three first properties satisfied by~$h$ by simply
calling~$h$ a 
{\em homotopy with image~$S$}, and the two last ones
by saying that~$h$ {\em preserves the functions belonging
to~$E$}
and {\em commutes with the elements of~$G$}.

\trois{comment-qproj}
The quasi-projectivity assumption can likely be removed, but it is currently needed in the proof 
for technical reasons. 

\trois{comment-equivariant}
The theorem does not
simply assert
the existence of a (model-theoretic)
deformation retraction of~$\hat V$ onto a tropical
sub-functor; 
it also ensures that this
deformation retraction can be required to
preserve finitely many arbitrary natural transformations 
from~$V$ to~$\Gamma_0$, and to commute with 
an arbitrary algebraic
action of a finite group. This strengthening 
of the expected statement is of course intrinsically
interesting, but this is not the only reason why Hrushovski 
and Loeser have decided to prove it.
Indeed, even if one
only wants to show the existence of a deformation retraction 
onto a tropical functor with no extra requirements,
there is a crucial step in the proof (by induction)
at which one needs to exhibit a deformation retraction
of a lower dimensional space to a tropical functor
{\em which preserves some natural transformations to~$\Gamma_0$
and the action of a suitable finite group}. 

\trois{annonce-preuve}
The remaining part of this section will now be devoted to 
(a sketch of) the proof of th.~\ref{main-theo}. We will first explain
(\ref{first-stepP1}-\ref{homo-curves})
what happens for curves; details can be found in~\cite{hl}, \S 7.

\deux{first-stepP1}{\bf The stable completion~$\widehat{\PP^1_k}$ is uniquely path-connected}. 

\trois{p1tree}
Let~$F\in \mathsf M$ and let~$x$ and~$y$ be two points
of~$\widehat{\PP^1_k}(F)$. One proves the following
statement, which one can roughly rephrase by saying that~$\widehat{\PP^1_k}$ is uniquely path-connected, or is a tree: 
{\em there exists a unique~$F$-definable sub-functor~$[x;y]$
of~$\widehat{\PP^1_F}$ homeomorphic
to a $k$-definable generalized segment with endpoints~$x$ and~$y$.} 

\trois{xy-explicit}
Let us describe~$[x;y]$
explicitly. 
If~$x=y$ there is nothing to do; if not, let
us distinguish two cases. In each of them, we will exhibit an~$F$-definable
generalized interval~$I$ and an~$F$-definable natural transformation
$\phi : I\to \widehat{\PP^1_F}$
inducing
a homeorphism between~$I$ and
a sub-functor of~$\widehat{\PP^1_F}$, and sending
 one of the endpoints of~$I$ to~$x$, 
and the other one to~$y$. 

\trois{cas-y-infty} {\em The case where, say, $y=\infty$}. One then has~$x=\eta_{a,r}$ for some~$a\in F$
and some~$r\in |F|$. We take for~$I$
the generalized interval~$[r;+\infty]$ defined in a natural way: 
if~$r>0$ it is homeomorphic to the interval~$[0;1/r]$; if $r=0$ one concatenates $[0;1]$ and $[1;+\infty]$. 
We take for~$\phi$
the natural transformation given by the formula~$s\mapsto \eta_{a,s}$, with the convention that~$\eta_{a,+\infty}=\infty$.

\medskip
\trois{cas-xy-finite}
{\em The case where~$x=\eta_{a,r}$ and~$y=\eta_{b,s}$ for~$a,b\in F$
and~$r,s\in |F|$}. We may assume that~$r\leq s$.

\begin{itemize}
\item[$\diamond$] If~$|a-b|\leq s$ then~$y=\eta_{a,s}$. We take for~$I$ the~$F$-definable interval~$[r; s]$
and for~$\phi$ the natural transformation given by the formula~$t\mapsto \eta_{a,t}$. 

\item[$\diamond$]  If~$|a-b|> s$ we take for~$I$ the concatenation of~$I_1:=[r;|a-b|]$ and of~$I_2=[|a-b|;s]$, 
and for~$\phi$ the natural transformation given by the formulas~$t\mapsto \eta_{a,t}$ for~$t\in I_1$ and~$t\mapsto \eta_{b,t}$
for~$t\in I_2$. 

\end{itemize}

\trois{convex-hull}
If~$\Delta$ is a finite,~$k$-definable (that is, Galois-invariant)
subset of~$\widehat{\PP^1_{\overline k}}$
the sub-functor~$\bigcup_{x,y\in \Delta}[x;y]$ 
of~$\widehat{\PP^1_k}$
is called the
{\em convex hull}
of~$\Delta$. It is a~$k$-definable ``twisted finite sub-tree"
of~$\widehat{\PP^1_k}$, where ``twisted" refers to the fact that the Galois action
on this sub-tree is possibly non-trivial. 

\deux{homo-p1} {\bf Retractions from~$\widehat{\PP^1_k}$ to twisted finite subtrees (\cite{hl}, \S 7.5).}

\trois{basic-formula-retraction}
Let~$U$ and~$V$ be
the~$k$-definable sub-functors of~$\PP^1_k$ respectively 
described by the conditions~$|T|\leq 1$ and~$|T|\geq 1$. Note that~$V$
is~$k$-definably
isomorphic to~$U$
through~$\psi : T\mapsto 1/T$; we still denote by~$\psi$ the induced
isomorphism~$\hat V\simeq \hat U$. 

\medskip
Let~$F\in \mathsf M$ and let~$x\in \hat U(F)$. One has~$x=\eta_{a,r}$
for
some~$a\in F\zero$ and some~$r\in |F\zero|$. For every~$t\in |F\zero|$, set
$$h(t,x)=\eta_{a, \max(t,r)},$$ and for every~$y\in \hat V(F)$, set
$$h(t,y)=\psi^{-1}(h(t,\psi(y)).$$ One immediately checks that both definitions
of~$h$ agree on~$\hat U\cap \hat V$, and we get by gluing 
a homotopy $h: [0;1]\times \widehat{\PP^1_k}\to  \widehat{\PP^1_k}$
with image~$\{\eta_{0,1}\}$. Hence~$\widehat{\PP^1_k}$ is ``$k$-definably
contractible".

\trois{retraction-subgraphs}
Let~$\Delta$ be a finite,~$k$-definable
subset of~$\widehat{\PP^1_{\overline k}}$, and let~$D$ be the convex hull
of~$\Delta\cup \{\eta_{0,1}\}$
(\ref{convex-hull}).
Define~$h_\Delta : [0;1]\times \widehat{\PP^1_k}\to  \widehat{\PP^1_k}$
as follows: for every~$x$, we denote by~$\tau_x$ the smallest time~$t$ such that~$h(t,x)\in D$, 
and we set~$h_\Delta(t,x)=h(t,x)$
if~$t\leq \tau_x$, and~$h_\Delta(t,x)=h(\tau_x,x)$ otherwise. The natural transformation~$h_\Delta$ is then a homotopy, 
whose image is the twisted polyhedron~$D$.

\deux{homo-curves}
{\bf Retractions from a curve to a twisted finite graph.}
In this paragraph we will sketch the proof of the theorem
when~$X=V$ is a
projective algebraic curve; the reader will find more details
in the survey \cite{bourb}, \S 4.2. 
One first chooses a finite $G$-equivariant
map~$f: X\to \PP^1_k$, inducing a
natural transformation~$\hat f : \hat X\to \widehat{\PP^1_k}$. 

\trois{key-point-curves}
The key point is the following (\cite{hl}, th. 7.5.1): {\em there exists a
finite~$k$-definable subset~$\Delta_0$ of~$\PP^1(\overline k)$ 
(or, in other words, a divisor on~$\PP^1_k$) such that for every
finite~$k$-definable subset~$\Delta$ of~$\PP^1(\overline k)$
containing~$\Delta_0$, the homotopy~$h_\Delta$ lifts {\em uniquely} to a homotopy
$h_\Delta^X: [0;1]\times \hat X\to \hat X$.} 
Hrushovski and Loeser prove it by carefully analyzing the behavior 
of the cardinality of the fibers $\widehat f$, as a function
from~$\widehat{\PP^1_k}$ to~$\NN$ (\cite{hl}, prop. 7.4.5). The definability of~$\hat X$ 
and~$\widehat{\PP^1_k}$ plays a crucial role for that purpose.

\trois{curves-conclusion}
Let us now explain why~\ref{key-point-curves}
allows to conclude. Let~$\Delta$ be as above, 
let~$D$ be the convex hull of~$\Delta\cup\{\eta_{0,1}\}$, and set~$D'=\hat f^{-1}(D)$. 

\medskip
$\bullet$ By a definability argument,
$D'$ is a ``twisted finite sub-graph" of~$\hat X$. 

\medskip
$\bullet$ By choosing~$\Delta$ sufficiently big, one may ensure that every function
belonging to~$E$
is locally constant outside~$D'$. This comes from the fact that every~$k$-definable
natural transformation from~$\hat X$ to~$\Gamma_0$ is locally constant outside a 
twisted finite subgraph of~$\hat X$: 
indeed, any~$k$-definable function can be described by
using only norms of regular functions; and the result for such a norm 
is deduced straightforwardly
from the behavior of~$|T|$ on~$\widehat{\PP^1_k}$, 
which is locally constant outside~$[0;\infty]$.  

\medskip
This implies that~$h^X_\Delta$ preserves the functions belonging to~$E$.
Moreover the uniqueness of 
$h^X_\Delta$ ensures that it 
commutes with the elements of~$G$, and we are done.

\trois{path-connected}
{\bf A consequence: path-connectedness of stable completions.}
Let~$W$ be a~$k$-definable sub-functor of a $k$-scheme of finite type, 
a $k$-definable sub-functor of $\Gamma_0^n$ or a $k$-generalized interval. 
We say that~$\hat W$
is {\em definably
path-connected}
if for every~$F\in \mathsf M$ and every~$(x,y)\in \hat W(F)^2$ 
there exists an~$F$-definable generalized interval~$I$ with endpoints~$o$
and~$e$, and an~$F$-definable continuous natural transformation~$h : I\to \hat W_F$
such that~$h(o)=x$ and~$h(e)=y$. 

\medskip
It follows easily
from the above that if~$W$ is a projective curve, then~$\hat W$
has
finitely many path-connected components, each of which is~$\bar k$-definable
(indeed, this holds for any twisted finite subgraph of $\hat W$). Starting from this result, Hrushovski and Loeser 
prove in fact that if~$Y$ is a~$k$-algebraic variety,~$\hat Y$ is definably path-connected
as soon as~$Y$ is geometrically connected
(\cite{hl}, th. 10.4.2); in general, there is a Galois-equivariant bijection between the set
of geometrical connected components of~$Y$ and that of path-connected components of~$\hat Y$.

\deux{cas-general}{\bf The general case of
th.~\ref{main-theo}: preliminaries.}

\trois{preli-maintheo-general} {\bf First reductions.}
By elementary geometrical arguments (\cite{hl}, beginning of \S 11.2)
one proves that there
exists a $G$-equivariant embedding of $X$ as a sub-scheme of an equidimensional projective $k$-variety $X^\sharp$
on which $G$ acts. Hence by replacing $X$ with $X^\sharp$ we may assume that
$X$ is projective and of pure dimension~$n$ for some~$n$; by extending the scalars to 
the perfect closure of~$k$ and by replacing~$X$ 
with its
underlying reduced sub-scheme, one can also assume that
its smooth locus is dense. {\em One proceeds then by induction on~$n$.}
The case~$n=0$ is obvious; now one assumes that~$n>0$ and that the theorem holds for smaller
integers.

\medskip
We will explain how to build a
homotopy from the {\em whole
space~$\widehat X$} to a tropical sub-functor
$S$ of $\widehat X$, 
which commutes with the elements of~$G$ and
preserves the characteristic function of~$V$ and the functions
belonging to~$E$; this homotopy will stabilize~$\hat V$ and the image
of $\hat V$ will be a $k$-definable sub-functor of $S$ (since $V$ is {\em strictly}
pro-definable), hence will be tropical as well. 

\medskip
By adding the characteristic function of $V$ to $E$
we thus reduce to the case where $V=X$.

\trois{blow-up}
{\bf Use of blowing-up to get a curve fibration}. 
One can blow up~$X$ along
{\em a finite set 
of closed points}
so that the resulting variety~$X'$ admits a morphism~$X'\to \PP^{n-1}_k$
whose generic fiber is a curve; {\em this is the point
where we need~$X$ to be projective}. By proceeding carefully
(and paying a special attention to the case where $k$ is finite) one can even ensure the following
(\cite{hl}, \S 11.2):

\medskip
$\bullet$ The action of~$G$ on~$X$ extends to~$X'$, and~$X'\to \PP^{n-1}_k$ is~$G$-equivariant. 

$\bullet$ There exists a $G$-equivariant
divisor~$D_0$ on~$X'$, finite over~$\PP^{n-1}_k$ and containing the exceptional divisor~$D$ of~$X'\to X$, 
and a~$G$-equivariant \'etale map
from~$X'\setminus D_0$ to~$\Aff^n_k$ (in particular, $D_0$ contains the singular locus of~$X'$). 

$\bullet$ There exists a non-empty open subset~$U$ of~$\PP^{n-1}_k$, whose pre-image
on~$X'$ is the complement of a divisor~$D_1$, 
and a factorization of~$X'\setminus D_1\to U$ through a finite, $G$-equivariant map~$f: X'\setminus D_1 \to \PP^1_k\times_k U$.

\medskip
If one builds a homotopy from~$\widehat{X'}$ to a twisted polyhedron
which commutes with
the action of~$G$, preserves the functions belonging to the (pullback of)~$E$, and the
characteristic function of
the exceptional
divisor $D$, it will descend to a homotopty 
on~$\widehat X$ satisfying the required properties, because every connected component of~$\widehat D$
collapses to a point.  Hence we reduce to the case where~$X'=X$.

\deux{concatenate}
{\bf The concatenation of a homotopy on the base and a fiberwise homotopy.}
This step is the core of the proof; we sum up here what is done in \cite{hl}, \S 11.3 and \S 11.4. 

\trois{base-fibre}
One first applies the {\em relative}
version of
the general construction we have described in~\ref{homo-curves}, to the finite map~$f$. It provides
a divisor~$\Delta$ on~$\PP^1_k \times_k U$
finite over~$U$
and a ``fiberwise'' homotopy~$h_\Delta$ 
on~$\widehat{\PP^1_k \times_k U}$ which lifts uniquely to a homotopy~$h_\Delta^X$ on~$\widehat{X\setminus D_1}$
whose image is a {\em relative (at most) one-dimensional tropical functor}
over~$\hat U$.  By taking~$\Delta$ big enough, one ensures that~$h_\Delta^X$ preserves the functions
belonging to~$E$ and commutes with the elements of~$G$. Moreover, we can also assume that the pre-image of~$\Delta$
on~$X\setminus D_1$ contains~$D_0\setminus D_1$. Under this last assumption, $\widehat{D_0\setminus D_1}$ is pointwise
fixed under~$h_\Delta^X$ at every time; therefore~$h_\Delta^X$ extends to a homotopy (which is still denoted by~$h_\Delta^X$) 
on~$\widehat{X\setminus D_1}\cup \widehat{D_0}$ which fixes~$\widehat D_0$ pointwise at every time. This homotopy
preserves the functions belonging to~$E$ and commutes with the action of~$G$, and its image~$\Upsilon$ is a
relative (at most)
one-dimensional tropical functor
over~$\widehat{\PP^{n-1}_k}$. Indeed, this is true by construction over~$\hat U$; and over the complement
of~$\hat U$, the fibers of
$\Upsilon$ coincide with those of~$\widehat{D_0}$, which are finite.

\trois{induction-hypothesis}
{\bf Use of the induction hypothesis}.
Our purpose is now to exhibit a homotopy~$h^\Upsilon$ on~$\Upsilon$, preserving the functions belonging to~$E$ and 
commuting with the action of~$G$, whose image is tropical 
of dimension at most~$n$.
The key point allowing such a construction 
is a theorem by Hrushovski and Loeser (\cite{hl}, th. 6.4.4)
which ensures that this problem is under algebraic control, in the following sense: 
there exists a finite quasi-Galois cover~$Z\to \PP^{n-1}_k$, with Galois group~$H$, and a finite family~$F$ of~$k$-definable natural transformation
from~$Z$ to~$\Gamma_0$, such that for every homotopy~$\lambda$ on~$\widehat{\PP^{n-1}_k}$ the following are equivalent: 

\medskip
i) $\lambda$ lifts to a homotopy on~$\hat Z$ preserving the functions belonging to~$F$; 

ii) $\lambda$ lifts to a homotopy on~$\Upsilon$ preserving the functions belonging to~$E$ and commuting with the action of~$G$. 

\medskip
Now {\em by induction hypothesis}
there exists a homotopy~$\lambda'$ on~$\hat Z$, preserving the functions belonging to~$F$ and commuting with the action of~$H$, 
and whose image is tropical of dimension~$\leq n-1$. Since~$\lambda'$ commutes with the action of~$H$, it descends
to a homotopy~$\lambda$ on~$\widehat{\PP^{n-1}_k}$. By its very definition, $\lambda$ satisfies condition~i) above, 
hence also condition~ii). Let~$h^\Upsilon$ be a homotopy on~$\Upsilon$ lifting~$\lambda$, preserving the functions belonging to~$E$ and commuting with the action of~$G$. 
Since the image of~$\lambda'$ is tropical 
of dimension~$\leq n-1$, so is the image of~$\lambda$ (because it is the quotient of that of~$\lambda'$ 
by the action of the finite group~$H$). As a consequence, the image~$\Sigma$ of~$h^\Upsilon$ is a relative (at most) one dimensional tropical
functor over a tropical functor of dimension~$\leq n-1$; therefore~$\Sigma$ is itself tropical of dimension~$\leq n$. 

\medskip
By making~$h^\Upsilon$ follow~$h_\Delta^X$, we get a homotopy~$h^0$ on~$\widehat {X\setminus D_1}\cup \widehat{D_0}$, 
preserving the functions belonging to~$E$ and commuting with the action of~$G$, whose image is the tropical functor~$\Sigma$.

\deux{inflation}
{\bf Fleeing away from~$\widehat {D_1}$.}

\trois{def-inflation}
{\bf The inflation homotopy: quick description}.
Hrushovski and Loeser define (\cite{hl}, lemma 10.3.2)
a ``inflation"
homotopy~$h^{\rm inf} : [0;1]\times \hat X \to \hat X$ which fixes
pointwise~$\widehat{D_0}$ at every time, preserves the functions
belonging to~$E$ and commutes with the action of~$G$,
and which is such that~$h^{\rm inf}(t,x)\in \hat X \setminus \widehat{D_1}$
for every~$x\notin \widehat{D_0}$ and every~$t>0$.  

\trois{const-inflation}
{\bf The inflation homotopy: construction}. 
One first defines a homotopy~$\alpha$ on~$\widehat{\Aff^n_k}$ by
``making the radii of balls increase'' (or, in other words, by generalizing to the higher dimension case
the formulas we have given for the affine line). It has the following genericity property: if~$x\in \widehat{\Aff^n_k}(F)$
for some~$F\in \mathsf M$, if~$Y$ is a~$(n-1)$-dimensional Zariski-closed subset of~$\Aff^n_F$ and
and if~$t$ is a non-zero element of~$|F\zero|$, then~$\alpha(t,x)\notin \hat Y(F)$. 
The homotopy~$\alpha$ is lifted to a homotopy~$\mu$
on~$\widehat{X\setminus D_0}$ thanks to the \'etale $G$-equivariant map~$X\setminus D_0 \to \Aff^n_k$, and the genericity 
property of~$\alpha$ is transferred to~$\mu$. 
The
homotpy~$h^{\rm inf}$ is then defined using a
suitable
``stopping time function''~$x\mapsto \tau_x$ from~$\widehat{X\setminus D_0}$
to~$\Gamma$ by the formulas: 

\medskip
$\bullet$ $h^{\rm inf}(t,x)=\mu(t,x)$ if~$x\notin \widehat{D_0}$ and if~$t\leq \tau_x$; 

$\bullet$ $h^{\rm inf}(t,x)=\mu(\tau_x,x)$ if~$x\notin \widehat{D_0}$ and if~$t\geq \tau_x$; 

$\bullet$ $h^{\rm inf}(t,x)=x$ if~$x\in \widehat{D_0}$.

\trois{verification-property-hinf}
Let us quickly explain why~$h^{\rm inf}$ satisfies the required properties.

\medskip
$\diamond$ Its continuity comes from the choice of the stopping time~$\tau_x$: the closer~$x$ is 
to~$\widehat{D_0}$, the smaller is~$\tau_x$. 

$\diamond$ The $G$-equivariance of~$h^{\rm inf}$ comes from its construction, and from
the~$G$-equivariance of~$X\setminus D_0 \to \Aff^n_k$. 

$\diamond$ The fact that $h^{\rm inf}(t,x)\in \hat X \setminus \widehat{D_1}$
for every~$x\notin \widehat{D_0}$ and every~$t>0$ is a particular case of the aforementioned
genericity property of~$\mu$.

$\diamond$ The fact that~$h^{\rm inf}$ preserves the functions belonging to~$E$
goes as follows. Those functions are defined using norms of regular functions; now if a regular function is invertible at a point~$x$
of~$\widehat{X\setminus D_0}$, 
then its norm is constant in a
neighborhood of~$x$, hence will be preserved by~$h^{\rm inf}(., x)$ if~$\tau_x$ is small enough; this is also obviously true
if the function vanishes in the neighborhood of~$x$. 

Of course, it may happen that the zero locus of a regular function involved in the description of~$E$ has some~$(n-1)$-dimensional irreducible component~$Y$,
and if~$x\in \hat Y\setminus \widehat{D_0}$, none of the above arguments will apply around
it. One overcomes this forthcoming
issue
by an additional work
at the very beginning of the proof: one
simply includes
such bad components
in the divisor~$D_0$. 

\trois{use-hinf}
By first applying~$h^{\rm inf}$, and then~$h^0$, one gets a continuous 
natural transformation~$h^1: J\times \hat X \to \hat X$ for some~$k$-definable generalized interval~$J$. Its image
is a~$k$-definable sub-functor~$\Sigma'$ of~$\Sigma$, hence is in
particular tropical of dimension~$\leq n$.

\deux{trop-homotopy}
{\bf The tropical homotopy and the end of the proof}. 

\trois{problem-h1}
The homotopy
$h^1$ of~\ref{use-hinf}
now enjoys all required properties, except (possibly) one of them: there is no reason
why~$\Sigma'$ should be pointwise fixed at every time, because~$h^{\rm inf}$ could disturb it. 

\trois{explanation-tropical}
Hrushovski and Loeser remedy it as follows. They proceed to the above construction
in such a way that there exists a $k$-definable $G$-equivariant sub-functor $\Sigma_0$ of $\Sigma$ satisfying the following conditions.

\medskip
1) {\em $\Sigma_0$ is pointwise fixed at every time under $h^1$}. This is achieved by ensuring that~$\Sigma_0$ is purely $n$-dimensional
and that there exists a finite extension $k'$ of $k$ and a $k'$-definable continuous
injection $\Sigma_0\hookrightarrow \Gamma_0^m$ 
that is preserved by $h^1$; this implies the required assertion by \cite{hl}, prop. 8.3.1. 

\medskip
2) There exists a homotopy  ~$h^{\rm trop}$ on~$\Sigma$, preserving the functions belonging to~$E$ and commuting with the action of~$G$, 
and whose image is precisely~$\Sigma_0$. The construction is purely tropical and rather technical, and we will not
give any detail here (see \cite{hl}, \S 11.5). 

\trois{conclusion-homotopy}
The expected homotopy~$h$ on~$\hat X$ is now defined by applying first~$h^1$, and then~$h^{\rm trop}$. 
Its image $S$ is equal to $h^{\rm trop}(\Sigma')$, hence is a $k$-definable sub-functor of $\Sigma_0$. By condition 1), it is
pointwise fixed at every time under $h^1$, hence under $h$.

%
%

\section{An application of the definability of~$\hat C$ for~$C$ a curve} 

\deux{intro-pl-geometry}
{\bf Presentation of the main result}. We fix a field $k$
which is complete with respect to a non-Archimedean absolute value. 

\trois{berko-domains}
Let us introduce some vocabulary (essentially coming 
from \cite{brk2}, \S 1.3). If $Z$ is a topological space
and if $Z'$ is a subset of $Z$, a family $(Z_i)$
of subsets of $Z'$ is said to be a {\em G-covering}
of $Z'$ (where G stands for ``Grothendieck") if every point $z$ of $Z'$ has a neighborhood 
in $Z'$ of the form $\bigcup_{i\in J}Z_i$ 
where $J$ is a {\em finite}
set of indices such that $z\in Z_i$ for every $i\in J$. 

\medskip
Any $k$-analytic space $X$ is equipped with a class of distinguished
compact subsets, the {\em affinoid domains}; the typical example
to have in mind is the subset of~$\Aff^{n,\rm an}_k$ defined by the inequalities
$$|T_1|\leq r_1\;{\rm and}\;\ldots\;{\rm and}\;|T_n|\leq r_n$$ where the $r_i$'s
are positive real numbers (this is the ``Berkovich compact polydisc of polyradius $(r_1,\ldots, r_n)$). 

\medskip
The subsets of $X$ that are G-covered
by the affinoid domains they contained are called {\em analytic domains}, and inherit
a canonical $k$-analytic structure. 
The open subsets of $X$ are analytic domains; so are the compact subsets of $X$ which are unions of finitely many
affinoid domains. 

\trois{pl-mult}
{\bf About piecewise-linear geometry (after Berkovich, \cite{sm2} \S1)}. 
We will follow Berkovich's convention in piecewise-linear 
geometry, which is based upon 
the {\em multiplicative notation},
in order to avoid using the use of
many somehow irrelevant
``$\log$" symbols. We will thus say
that a subset of~$(\RR_+^\times)^n$ is a polyhedron
if its image under $\log$ is a polyhedron
of~$\RR^n$ (def.~\ref{polyhedron}), and that a map between
two polyhedra is piecewise-linear if it is piecewise $\QQ$-linear 
in the usual sense modulo
$\log$ isomorphisms. 

A {\em polyhedral structure} on a compact
topological space $P$ is a map $P\to (\RR^\times_+)^n$ 
(for some $n$) inducing an homeomorphism between $P$ and a polyhedron 
of $(\RR^\times_+)^n$
Two polyhedral 
structure $i\colon P\to (\RR^\times_+)^n$ and $j\colon P\to (\RR^\times_+)^m$ 
are said to be {\em equivalent}
if there exists a PL homeomorphism $u\colon i(P)\simeq j(P)$ such that $j=u\circ i$. 

A {\em polyhedral chart} on a topological space $Z$ is a compact
subset $P$ of $Z$ (the {\em support} of the chart)
equipped with a class of polyhedral structures. 
Two polyhedral 
charts $i\colon P\simeq P_0$ and $j\colon Q\simeq Q_0$ are {\em compatible}
if both $i(P\cap Q)$ and $j(P\cap Q)$ are polyhedra (this depends only of the classes of
$i$ and $j$). 

A {\em polyhedral atlas}
on a topological space $S$ is a family of
polyhedral charts which are pairwise compatible and whose support G-cover $S$.
Two polyhedral atlases on $S$ are said
to be {\em equivalent} if their union is a polyhedral atlas, that is, 
if any chart of the first one is compatible with any chart of the second one. A {\em piecewise-linear space}, 
or PL space for short, is a topological space equipped with an equivalence class of polyhedral atlases, 
or with a maximal polyhedral atlas -- this amounts to the same. When we will speak about 
a polyhedral atlas on a given PL space, we will implicitly assume that it belongs to the equivalence
class defining the PL structure involved.

If $S$ is a PL space a {\em PL subspace}
of $S$ is a subset of $S$ which is G-covered by (support of) charts of the maximal
polyhedral atlas of $S$ ; such a subspace $\Sigma$ inherits a canonical PL structure. Open subsets of $S$
and compact subsets which are the union of finitely many (support of) charts of the maximal polyhedral atlas
of $S$ are examples of PL subspaces of $S$. 

A map $f\colon S\to\Sigma$ bewteen PL spaces is said to be PL if there is a polyhedral atlas $\mathfrak A$ on $S$
such that for every $P\in \mathfrak A$ the image $f(P)$ belongs to the maximal polyhedral
atlas of $\Sigma$, 
and $P\to f(P)$ is PL modulo the identifications of its source and its target with polyhedra. If moreover
$\mathfrak A$ can be chosen
so that $P\to f(P)$ is injective for every $P\in \mathfrak A$, we say that $f$ is a {\em piecewise immersion}. 

\trois{def-sn}
Let $n\in \NN$. Once equipped with the family of all its polyhedra, 
$(\RR^\times_+)^n$ becomes a 
a PL space. 

Let us denote by~$S_n$ the ``skeleton"
of~${\mathbb G}_m^{n,\rm an}$. This is the set of semi-norms of the form
$\eta_{r_1,\ldots, r_n}:=\sum a_{I}T^I\mapsto \max |a_I| r^I$. 
The very definition of~$S_n$ provides it with a homeomorphism
onto~$(\RR^\times_+)^n$, through which it inherits a
piecewise-linear structure. 

\trois{imrec-skel-state}
Let~$X$ be a~$k$-analytic space of dimension~$n$, and let~$f: X\to {\mathbb G}_m^{n,\rm an}$
be a morphism. In~\cite{polyt} (see also the {\em erratum}~\cite{errat-polyt})
the author has proven that~$f^{-1}(S_n)$ inherits a canonical
PL structure, with respect to which $f^{-1}(S_n)\to S_n$
is a piecewise immersion. This generalizes his preceding work~\cite{squel}, 
in which some additional assumptions on~$X$ were needed, and in which the canonicity
of the PL structure (answering a question by Temkin) had not been addressed. 

\medskip
Moreover, the proof given in~\cite{squel} used de Jong's alterations, while that 
given in~\cite{polyt} replaces it by the $k$-definability of~$\hat C$ for~$C$ an algebraic curve,
(th.~\ref{hatv-prodef}) and which simply comes from Riemann-Roch theorem for curves --,
together with some general results 
about model theory of valued fields that are established 
by Haskell, Hrushovski and Macpherson
in \cite{hhm}. 

\medskip
The purpose of
this section is to give a very rough sketch of the latter proof. 

\deux{czn-pl-berko}
{\bf Canonical PL subsets of Berkovich spaces.} 
Let~$X$ be a Hausdorff~$k$-analytic space of dimension~$n$.

\trois{abhyankar-point}
Let $x\in X$, let $d$ be transcendance degree
of $\widetilde{\mathscr H(x)}$ over $\widetilde k$
and let $r$ be the rational rank of the abelian group 
$|\mathscr H(x)^\times|/|k^\times|$. One alway has
the inequality $d+r\leq n$; the point $x$ will be said to be {\em Abhyankar}
if $d+r=n$.

\trois{def-skeleton}
Let~$S$ be a locally closed subset of~$X$. 
We will say that~$S$ is a {\em skeleton}
if it consists of Abhyankar points 
and if it 
it admits a PL structure satisfying the following conditions. 

\medskip
\begin{itemize}

\item[$\bullet$] For every analytic domain~$Y$ of~$X$ and every invertible function~$f$ on~$Y$, the
intersection~$Y\cap S$ is a PL subspace of~$S$, and the restriction of~$|f|$
to~$Y\cap S$ is PL. 

\medskip
\item[$\bullet$] There exists a polyhedral atlas $\mathfrak A$ 
on $S$ having the following property: for every $P\in \mathfrak A$, 
the class of polyhedral structures defining the chart $P$
admits a representative of the form $(|f_1|, \ldots, |f_m|)|_P$ where
the $f_i$'s are invertible functions on an analytic domain $Y$ of $X$ which contains $P$.

\end{itemize}

\medskip
If such a structure exists on~$S$, it is easily seen to be unique. The archetypal
example of a skeleton is $S_n$ (\ref{def-sn}).

\trois{criterion-skeleton}
{\bf A simple criterion for being a skeleton.}
The conditions for being a skeleton may seem slightly complicated to check. But in fact, 
in practice
this is not so difficult. Indeed, let
us assume that~$X$ is affinoid and irreducible, 
of dimension~$n$, and let~$P$ be a compact subset of~$X$ consisting 
of Abhyankar points. The latter property ensures that every point of~$P$ is Zariski-generic. 
In particular, any non-zero analytic function on~$X$ is invertible on~$P$.

Now assume that for every finite family~$(f_1,\ldots, f_m)$ of non-zero analytic functions on~$X$,
the image~$(|f_1|, \ldots, |f_m|)(P)$ is a polyhedron of~$(\RR^\times_+)^m$, and that there exists such a family
$(f_1,\ldots, f_m)$ with~$(|f_1|, \ldots, |f_m|)|_P$ injective. Under this assumption, $P$ is a skeleton
(\cite{polyt}, Lemme 4.4). 

\medskip
The two
main
ingredients of the proofs are the Gerritzen-Grauert theorem
(that describes the affinoid domains of $X$, see for instance
\cite{semialg}, lemme 2.4), and the density of rational functions on~$X$ inside the ring
of analytic functions of any rational affinoid domain of~$X$. 

\trois{rephrase}
We can now rephrase th.~\ref{imrec-skel-state}
as follows. 

\trois{thm-rephrased}
{\bf Theorem (\cite{polyt}, th. 5.1)}. 
{\em If
$f$ 
is any morphism from $X$ to ${\mathbb G}_m^{n,\rm an}$ then $f^{-1}(S_n)$ is a skeleton,
and 
$f^{-1}(S_n)\to S_n$ is a piecewise immersion}.

\deux{gen-strategy-pl}
{\bf The general strategy: algebraization and choice of a faithful tropicalization.}

\trois{algebraization}{\bf Algebraization.} 
In order to prove th.~\ref{thm-rephrased}
we first algebraize the situation
by standard arguments which we now outline
(for more details, see \cite{polyt} \S 5.4 and \S 5.5). 

\medskip
First of all, one extends the scalars to the
perfect closure of~$k$, and then replaces~$X$ with its underlying reduced space; 
we can do it because the expected assertions are insensitive to radicial extensions and 
nilpotents. Now~$X$ is generically quasi-smooth 
(quasi-smoothness is the Berkovich
version of rig-smoothness, see \cite{flat}, \S 4.2 {\em et sq.}). 
Every point of $f^{-1}(S_n)$ is Abhyankar (because so are the points
of $S_n$ and because $\dim X=n$), hence is quasi-smooth; one can thus
shrink~$X$ so that it is 
itself quasi-smooth. Krasner's lemma
then ensures that~$X$ is G-locally algebraizable
(see for instance \cite{polyt}, (0.21))
which eventually allows, 
because being a skeleton is easily seen to be a G-local
property,
to reduce
to the case where~$X$
is a connected, irreducible
rational 
affinoid domain
of~$\sch X^{\rm an}$ for~$\sch X$ an irreducible, normal
(and even smooth)
algebraic variety over~$k$. 

\medskip
By density arguments,
we also can assume
that~$f$ is induced by a dominant, generically finite
algebraic map
from~$\sch X$
to
${\mathbb G}_m^n$
(which we still denote by~$f$). Now let~$x\in f^{-1}(S_n)$. It is Abhyankar
hence Zariski-generic. 
By openness of finite, flat morphisms, 
it therefore admits a connected
affinoid
 neighborhood~$U_x$ in~$\sch X^{\rm an}$ 
such that~$f$ 
induces a finite and flat map from~$U_x$ to an affinoid domain~$V_x$ of~${\mathbb G}_m^{n,\rm an}$ 
with~$V_x\cap S_n$ being a non-empty $n$-dimensional simplex. 
Now there are finitely many such~$U_x$'s covering~$X\cap f^{-1}(S_n)$. 
It is then sufficient to prove that for every~$x\in S_n$
the intersection $U_x\cap f^{-1}(S_n)$
is a skeleton and that $U_x\cap f^{-1}(S_n)\to S_n$ is a piecewise-immersion.
Therefore we may assume that~$f$ 
induces a finite and flat map from~$X$ to an affinoid domain~$V$ of~${\mathbb G}_m^{n,\rm an}$ 
with~$V\cap S_n$ being a non-empty $n$-dimensional simplex. 

\medskip
The key result is now the following
theorem. 

\trois{fin-separation-thm}
{\bf Finite separation theorem.}
{\em There exist finitely many non-zero
rational functions~$g_1,\ldots, g_r$
on~$\sch X$ whose norms separate the pre-images of~$x$
under~$f$
for every~$x\in S_n$.} 

\medskip
We postpone the outline of the proof of th.~\ref{fin-separation-thm} to paragraph~\ref{separgauss}; 
we are first going to explain how one can use it to prove th.~\ref{imrec-skel-state}, 
or more precisely its rephrasing written down at~\ref{rephrase}.

\trois{trop-dim}
{\bf Pre-image of~$S_n$ and tropical dimension}
Let~$f_1,\ldots, f_n$ be the invertible functions on~$\sch X$ that define~$f$. 
For every compact analytic domain $Y$ of $\sch X\an$, the subset 
~$(|f_1|, \ldots, |f_n|)(Y)$ is a polyhedron ({\em cf.} \cite{sm2}, cor. 6.2.2; see also
\cite{polyt}, th. 3.2). 

If~$x\in \sch X^{\rm an}$, the {\em tropical dimension}
of~$f$ at~$x$ is the infimum of the dimensions of the polyhedra~$(|f_1|, \ldots, |f_n|)(Y)$ 
for~$Y$ going through the set of compact analytic neighborhoods of~$x$. One can then characterize~$f^{-1}(S_n)$ as the
subset
of~$\sch X^{\rm an}$
consisting
of
points at which the tropical dimension of~$f$ is exactly~$n$
(this follows from \cite{polyt}, th. 3.4). 

\trois{intersection-skeleton}
{\bf The compact
subset $f^{-1}(S_n)\cap X$ of $\sch X\an$ is a skeleton, and $U_x\cap f^{-1}(S_n)\to S_n$ is a piecewise-immersion.} 
Since every point of~$f^{-1}(S_n)$
is Zariski-generic, the functions~$g_i$'s are invertible on~$f^{-1}(S_n)$; hence we can shrink
$X$ so that
the~$g_i$'s are invertible on
it. Let~$h_1,\ldots, h_m$ be arbitrary non-zero analytic functions on~$X$. 

We set for short $|f|=(|f_1|,\ldots, |f_n|) \colon X \to (\RR^\times_+)^n$
and define $|g|$ and $|h|$ analogously. 
Let $\pi$
be the map $(|f|,|g|,|h|)$ from $X$ to $(\RR^\times_+)^{n+m+r}$; the
image $\pi(X)$ is a polyhedron. 

The restriction of $|f|$ to $S_n$ is injective
by the very definition of $S_n$, and $|g|$ separates the pre-images of every point
of $S_n$ on $\sch X\an$; therefore the restriction of $(|f|,|g|)$ to $f^{-1}(S_n)$
is injective; hence~$\pi|_{f^{-1}(S_n)\cap X}$ is injective. 

Let us choose a triangulation of
 the compact polyhedron~$\pi(X)$ by convex compact 
 polyhedra, and let~$P$ be the union of the~$n$-dimensional closed
 cells~$Q$ such that the following holds: {\em the restriction of~$(|f_1|, \ldots, |f_n|)$ to~$Q$ is injective}.
 Using the characterization of~$f^{-1}(S_n)$
 through tropical dimension (\ref{trop-dim}), the openness of finite, flat morphism 
 and the fact that~$V\cap S_n$ is non-empty of pure dimension~$n$, one proves that~$\pi(f^{-1}(S_n)\cap X)=P$
 (\cite{polyt}, 5.5.1); in particular
 this is a polyhedron and in view of~\ref{criterion-skeleton}
this implies that $f^{-1}(S_n)\cap X$ is a skeleton. Moreover
modulo the isomorphism $(f^{-1}(S_n)\cap X)\simeq P$  induced by $\pi$, the function $|f|$
is nothing but the projection to the first $n$ variables; it is PL and injective on each cell of $P$, 
which implies
that $f^{-1}(S_n)\cap X\to S_n$ 
is a piecewise immersion. This ends the proof
of th.~\ref{thm-rephrased}. 

\deux{separgauss}
{\bf About the proof of the finite separation theorem \ref{fin-separation-thm}.}

\trois{rapp-separ-th}
This theorem
asserts that there exist
finitely many
elements~$g_1,\ldots, g_r$ of the function field~$k(\sch X)$
which separate the extensions of every real-valued Gau\ss norm 
on~$k(T_1,\ldots, T_n)$. 
In fact, we establish the more general, purely valuation-theoretic  following theorem
(\cite{polyt}, th. 2.8). Let~$k$ be an
{\em arbitrary} valued field, let~$n$
be an integer, and
let~$L$ be a finite extension of~$k(T_1,\ldots, T_n)$. There exists a finite subset~$E$ of~$L$ such that the following hold: {\em for every ordered
abelian group~$G$
containing~$|k|$ and any~$n$-uple~$r=(r_1,\ldots, r_n)$ of elements of~$G$, the elements of~$E$ separate the extensions of the $G$-valued
valuation~$\eta_r$
of~$k(T_1,\ldots T_n)$ to~$L$.}

\trois{part-case}
{\bf A particular case}. 
In order to illustrate the general idea the proof
is based upon, 
let us explain
what is going on
in a simple particular case,
namely if~$k$ is algebraically closed, and if~$n=1$; we write~$T$
instead of~$T_1$. The
field~$L$ can then be written~$k(C)$ for a suitable projective, smooth, irreducible curve~$C$
over~$k$, equipped with a finite
morphism~$C\to \PP^1_k$ inducing the extension~$k(T)\hookrightarrow L$. 
The map~$C\to \PP^1_k$ induces a~$k$-definable natural transformation
 $\hat f : \hat C \to \widehat{\PP^1_k}$ at the level of stable
completions; it follows
from th.~\ref{hatv-prodef} 
that $\hat C$ and $\widehat{\PP^1_k}$ are $k$-definable
(for $\widehat{\PP^1_k}$ this can also be seen directly, {\em cf.}~\ref{example-1}).

\medskip
The map~$r\mapsto \eta_r$ defines a~$k$-definable natural embedding
$\Gamma\hookrightarrow \widehat{\PP^1_k}$; let~$\Delta$ be the pre-image of~$\Gamma$
under~$\hat f$.This is a sub-functor of~$\hat C$ which is~$k$-definable by its very definition, and 
the fibers of~$\Delta\to \Gamma$ are finite. 

\medskip
Let~$F$ be an algebraically closed valued extension of~$k$, and let~$r\in |F|$. The pre-image~$D$
of~$\eta_{r,F}$ inside~$\Delta(F)$ is by construction definable over the set of parameters~$k\cup \{r\}$, and is finite. 
By a result
proven in~\cite{hhm},
this implies that every element of~$D$ is 
{\em individually}
$(k\cup\{r\})$-definable\footnote{It is crucial for this result that the additional parameter~$r$ belongs to the value group
of~$F$, and not to~$F$ itself. For instance, let~$a$ be an element of~$F$ which is transcendent over~$k$, assume that char.~$k\neq 2$ and let~$E$ be the two-element 
set of square roots of~$a$. Then~$E$ is globally definable over~$k\cup\{a\}$, but this is not
the case of any of those two square roots: one cannot distinguish
between them by a formula involving only
$a$ and elements of~$k$.}. In other words, for every~$d\in D$, 
there exists a~$k$-definable natural transformation~$\sigma : \Gamma\to \Delta$ such that~$d=\sigma(\eta_{r,F})$. 

\medskip
As this holds for arbitrary~$F$ and~$r$, the celebrated compactness theorem of model theory
ensures the existence of finitely many sections~$\sigma _1, \ldots, \sigma_m$ 
of~$\Delta\to \Gamma$ such that~$\Delta=\bigcup \sigma_i(\Gamma)$. Every sub-functor~$\sigma_i(\Gamma)$
is~$k$-definable, and~$k$-definably isomorphic to~$\Gamma$ through~$\sigma_i$. 
One then
easily builds, starting from the equality~$\Delta=\bigcup \sigma_i(\Gamma)$, a~$k$-definable
embedding $\Delta\hookrightarrow \Gamma^N$ 
(for some big enough~$N$ -- one simply has to be able to realize the ``coincidence diagram" of
the~$\sigma_i$'s inside~$\Gamma^N$). 

\medskip
Now Hrushovski and Loeser have proven that is~$S$ is a~$k$-definable
sub-functor of~$\hat C$
such that there exists a~$k$-definable
embedding $S\hookrightarrow \Gamma^N$ then
then there exists such an isomorphism
{\em induced by
the norms of finitely many~$k$-rational functions on~$\hat C$}
(\cite{hl}, prop. 6.2.7);
this comes from the explicit description of~$\hat C$ as a~$k$-definable functor, and from general
results established in~\cite{hhm}. 
Applying this to~$\Delta$, we get the existence
of finitely many rational functions
on~$C$ whose norms
separate points of~$\Delta$; in particular those functions
separate for every algebraically closed
valued extension~$F$ of~$k$ and every~$r\in |F^*|$ the pre-images of~$\eta_{r,F}$ in~$\hat C(F)$. They
{\em a fortiori}
separate the extensions of the valuation~$\eta_{r,k}$ to the field~$L$
because  
by elementary valuation-theoretic arguments
such an extension is 
always the restriction of a valuation on~$F(C_F)=L\otimes_{k(T)}F(T)$ inducing~$\eta_{r, F}$ on~$F(T)$.

\trois{separate-fin-general}
{\bf The general case.} 
The proof
 goes by induction on~$n$. The crucial step is of course
the one that consists in going from~$n-1$
to~$n$, and it roughly
consists of
a relative version of what we have
done above.

\deux{complements-pl}
{\bf Some complements.}
We still  denote by $k$ a complete, non-Archimedean field and by ~$X$
an~$n$-dimensional~$k$-analytic space.

\trois{finite-union-premi}
{\bf Finite union of pre-images of the skeleton.} Let~$f_1,\ldots, f_m$ be morphisms
from~$X$
to~${\mathbb G}^{n,\rm an}_m$. 

\medskip
Th.~\ref{thm-rephrased}
states that~$f_i^{-1}(S_n)$ is a skeleton for every~$i$. Th. 5.1 of~\cite{polyt}
in fact also states that the union~$\bigcup_i f_i^{-1}(S_n)$ is still a skeleton; this essentially means that for every~$(i,j)$, the intersection
of the skeletons $f_i^{-1}(S_n)$
and~$f_j^{-1}(S_n)$
is PL in both of them. The proof consists in exhibiting in a rather explicit way an integer~$N$, 
a ``Shilov section with non-constant radius" $\sigma : X\to \Aff^N_X$, 
and a skeleton~$P\subset \Aff^N_X$ (described as the 
pre-image of $S_{N+n}$ under a suitable map)
such that
for every~$i$, the section~$\sigma$
identifies~$f_i^{-1}(S_n)$ with a PL subspace of~$P$
(\cite{polyt}, 5.6.2). 

\trois{stabilize}
{\bf Stabilization after a finite, separable extension.} 
For every complete extension~$F$ of~$k$, let~$\Sigma_F\subset X_F$ be
the skeleton~$\bigcup_i f_{i,F}^{-1}(S_{n,F})$, where~$S_{n,F}$ is the skeleton of~${\mathbb G}^{n,\rm an}_{m,F}$. 

If~$F\hookrightarrow L$, there is a natural surjection~$\Sigma_L\to \Sigma_F$, which is a piecewise immersion
of PL spaces. Th. 5.1 of~\cite{polyt} in fact also states
that there exists a finite separable extension~$F$ of~$k$ such that~$\Sigma_L\to \Sigma_F$
is a homeomorphism for every complete extension~$L$ of~$k$. To see it, one essentially refines the proof of the ``separation theorem" 
on Gau\ss~
norms that is sketched at~\ref{separgauss}
{\em et sq.}
to show that after making a finite, separable extension of the ground field, one can exhibit a finite family of functions
separating {\em universally}
(that is, after any extension of the ground field) the extensions of Gau\ss~
norms; this is part of aforementioned th. 2.8 of {\em loc. cit}.

\end{document}